\documentclass[11pt,twoside]{article}  

\usepackage[margin=1.25in]{geometry}

\usepackage{amsmath,amssymb,amsfonts,amsthm}
\usepackage{comment}
\usepackage[numbers]{natbib}
\usepackage{xspace}
\usepackage{graphicx}
\usepackage{color}
\usepackage{appendix}
\usepackage{hyperref}

\newtheorem{lem}{Lemma}
\newtheorem{thm}{Theorem}
\newtheorem{cor}{Corollary}

\newtheorem{defenv}{Definition}

\newcommand{\nulls}{\mathcal{H}_{0}}
\newcommand{\sigs}{\mathcal{H}_{1}}

\newcommand{\E}{\mathbb{E}}

\newcommand{\Var}{\mathrm{Var}}
\renewcommand{\P}{\mathbb{P}}
\newcommand{\R}{\mathbb{R}}

\newcommand{\N}{\mathbb{N}}
\newcommand{\Norm}{\mathcal{N}}
\newcommand{\FDR}{\mathrm{FDR}}

\newcommand{\FNR}{\mathrm{FNR}}

\newcommand{\FDP}{\mathrm{FDP}}
\newcommand{\FNP}{\mathrm{FNP}}

\newcommand{\discoveries}{\mathcal{I}}
\newcommand{\numobs}{n}
\newcommand{\card}[1]{\ensuremath{\mbox{card}(#1)}}
\newcommand{\Xseq}{\ensuremath{X_1^\numobs}}

\newcommand{\event}{\mathcal{E}}
\newcommand{\eventwrap}[1]{\left\lbrace{#1}\right\rbrace}
\newcommand{\kstar}{k^{\ast}}
\newcommand{\kstarp}{k_{+}^{\ast}}
\newcommand{\kstarm}{k_{-}^{\ast}}
\newcommand{\lstar}{\ell^{\ast}}
\newcommand{\lstarp}{\ell_{+}^{\ast}}
\newcommand{\lstarm}{\ell_{-}^{\ast}}

\newcommand{\FDPm}{\FDP_{-}^{\ast}}
\newcommand{\FNPm}{\FNP_{-}^{\ast}}
\newcommand{\FDPp}{\FDP_{+}^{\ast}}
\newcommand{\FNPp}{\FNP_{+}^{\ast}}

\newcommand{\rhonull}{\rho_{0}}
\newcommand{\rhosig}{\rho_{1}}
\newcommand{\rhocross}{\rho_{\mathrm{c}}}
\newcommand{\Vcross}{V_{\mathrm{c}}}
\newcommand{\vcross}{v_{\mathrm{c}}}

\newcommand{\ordD}{\Delta}
\newcommand{\transf}{f}

\newcommand{\epsinv}{\epsilon^{-1}}
\newcommand{\cproxy}{c_{0}}
\newcommand{\PW}{\P_{n}}
\newcommand{\Topk}{Top-$K$\xspace}
\newcommand{\topk}{top-$K$\xspace}
\newcommand{\nullDelta}{\Delta_{0}}
\newcommand{\sigDelta}{\Delta_{1}}
\newcommand{\model}{\mathbb{M}}
\newcommand{\allmodels}{\mathcal{M}}

\newcommand{\falsedisc}[1]{L\left({#1}\right)}
\newcommand{\bandevent}{\event_{\mathrm{band}}}
\newcommand{\proxyeventm}{\event_{\mathrm{proxy},-}}
\newcommand{\proxyeventp}{\event_{\mathrm{proxy},+}}
\newcommand{\goodeventm}{\event_{-}}
\newcommand{\goodeventp}{\event_{+}}

\newcommand{\gaussord}{\chi}

\newcommand{\Qs}{\mathcal{Q}}
\newcommand{\BH}{\mathrm{BH}}
\newcommand{\lo}{\mathrm{lo}}

\newcommand{\widgraph}[2]{\includegraphics[keepaspectratio,width=#1]{#2}}

\newcommand{\mydefn}{\ensuremath{: \, =}}
\newcommand{\defn}{\mydefn}

\newcommand{\fdr}{\ensuremath{\alpha}}
\newcommand{\fnr}{\ensuremath{\beta}}

\newcommand{\lehmann}{\ensuremath{\gamma}}

\newcommand{\betapar}{\ensuremath{s}}
\newcommand{\etapar}{\ensuremath{t}}
\newcommand{\EPSINV}{\ensuremath{\frac{1}{\epsilon}}}
\newcommand{\mprob}{\ensuremath{\mathbb{P}}}

\long\def\comment#1{}

\makeatletter
\long\def\@makecaption#1#2{
        \vskip 0.8ex
        \setbox\@tempboxa\hbox{\small {\bf #1:} #2}
        \parindent 1.5em  
        \dimen0=\hsize
        \advance\dimen0 by -3em
        \ifdim \wd\@tempboxa >\dimen0
                \hbox to \hsize{
                        \parindent 0em
                        \hfil 
                        \parbox{\dimen0}{\def\baselinestretch{0.96}\small
                                {\bf #1.} #2
                                } 
                        \hfil}
        \else \hbox to \hsize{\hfil \box\@tempboxa \hfil}
        \fi
        }
\makeatother

\begin{document}

\begin{center}

  {\bf{\LARGE{Lower bounds in multiple testing: A framework based on derandomized
        proxies}}}
  
\vspace*{.2in}

{\large{
\begin{tabular}{c}
Maxim Rabinovich$^\dagger$, Michael
I. Jordan$^{*, \dagger}$, Martin J. Wainwright$^{*, \dagger, \ddagger}$ 
\\
\\
\texttt{\{rabinovich,jordan,wainwrig\}@berkeley.edu}\\
\\
\end{tabular}
\begin{tabular}{ccc}
  Departments of Statistics$^*$ and EECS$^\dagger$, University of California, Berkeley \\
  Voleon Group$^\ddagger$, Berkeley
\end{tabular}
}}

\vspace*{.2in} \today
\vspace*{.2in}

\begin{abstract}
The large bulk of work in multiple testing has focused on specifying
procedures that control the false discovery rate (FDR), with
relatively less attention being paid to the corresponding Type II
error known as the false non-discovery rate (FNR).  A line of more
recent work in multiple testing has begun to investigate the tradeoffs
between the FDR and FNR and to provide lower bounds on the performance
of procedures that depend on the model structure.  Lacking thus far,
however, has been a general approach to obtaining lower bounds for a
broad class of models.  This paper introduces an analysis strategy
based on derandomization, illustrated by applications to various
concrete models.  Our main result is meta-theorem that gives a general
recipe for obtaining lower bounds on the combination of FDR and FNR.
We illustrate this meta-theorem by deriving explicit bounds for
several models, including instances with dependence, scale-transformed
alternatives, and non-Gaussian-like distributions. We provide
numerical simulations of some of these lower bounds, and show a close
relation to the actual performance of the Benjamini-Hochberg (BH)
algorithm.

\end{abstract}
\end{center}


\section{Introduction}

The past decades have witnessed a tremendous amount of research on
control of the false discovery rate (FDR), an analogue of type I error
in multiple testing problems (e.g,~\cite{Ben95FDR, BH97,
  BY01,foster2008alpha, Storey02, Storey04, Gen06FDRWeights,
  genovese2002operating}).  The large majority of work has focused on
developing procedures that are guaranteed to control the FDR at a
pre-specified level, under various assumptions on the structure of the
$p$-values.  The literature on classical hypothesis testing is replete
with combined analyses of the type I error and statistical power for
specific classes of models.  More recently, analogues of such analyses
have begun to appear in the FDR literature.  Using the false
non-discovery rate (FNR)---the fraction of tests in which the null is
incorrectly not rejected---as a measure of the (lack of) power,
several authors have established lower bounds on combinations of the
FDR and FNR~\citep{Ari16DistFreeFDR,Ji12UPS,Ji14Optimal}, and most
recently, a non-asymptotic tradeoff between FDR and FNR has been
established in the same setting~\citep{Rab20Tradeoff}.

Most results of this type occur within models where test statistics
arise from Gaussian-like location models with independent
observations, meaning that the alternatives are assumed to be location
shifts of the null, and the noise variables are independent with
Gaussian-type tail behavior. Models of this type have the advantage of
analytical tractability, while still permitting the expression of
central features of many multiple testing problems. Many results apply
to variants of the rare-weak (RW) model, in which problem difficulty
is parameterized by the \emph{rarity} of signals and their
\emph{weakness}.  This model was initially introduced for studies
multiple testing using the family-wise error rate (FWER) as the Type-I
error concept~\citep{Don04HC,Don15HC,Jin14RareWeak},

Unfortunately, the analysis techniques underlying existing results do
not extend much beyond the setting of independent and Gaussian-like
test statistics. This limitation appears most dramatically when one
seeks to derive non-asymptotic guarantees, as these results depend
critically on the ability to control tail probabilities and apply
concentration inequalities.  In this work, we introduce a modeling and
analysis strategy for lower bounds on the FDR and FNR that is far less
dependent on independence and analytical tractability. Building on the
proof strategy introduced in our earlier paper~\citep{Rab20Tradeoff},
we show how it is fruitful to perform a version of
\emph{derandomization}---that is, we relate a given multiple testing
procedure to a derandomized version that always makes a \emph{fixed}
number of discoveries. Given a proxy $\kstar$ for the number of
discoveries, we show how to further eliminate the randomness
associated with the location of the data to obtain deterministic
proxies $\lstar$ for lower bounding the number of \emph{false}
discoveries with constant probability. Using the two proxies $\kstar$
and $\lstar$ together yields constant-probability lower bounds on the
proportion of false discoveries and false non-discoveries, and these
translate directly into corresponding lower bounds on the FDR and FNR
that hold for any procedure applied to the given model.

With this context, the central contribution of this paper is a
meta-theorem that establishes a tradeoff between FDR and FNR in terms
of these proxies for a broad class of models defined in
Section~\ref{SecModel} In contrast to minimax theory for estimation,
where tools like Le Cam's method, Fano's method, and Assouad's lemma
(e.g., see~\cite{Yu97,Tsy08Estimation, Wai19} for background) play the
role of meta-theorems that can be instantiated, such general tools are
currently absent from the multiple testing literature. The main
theorem of this paper appears to be the first of its kind for
controlling FDR and FNR in multiple testing.

\begin{figure}[t]
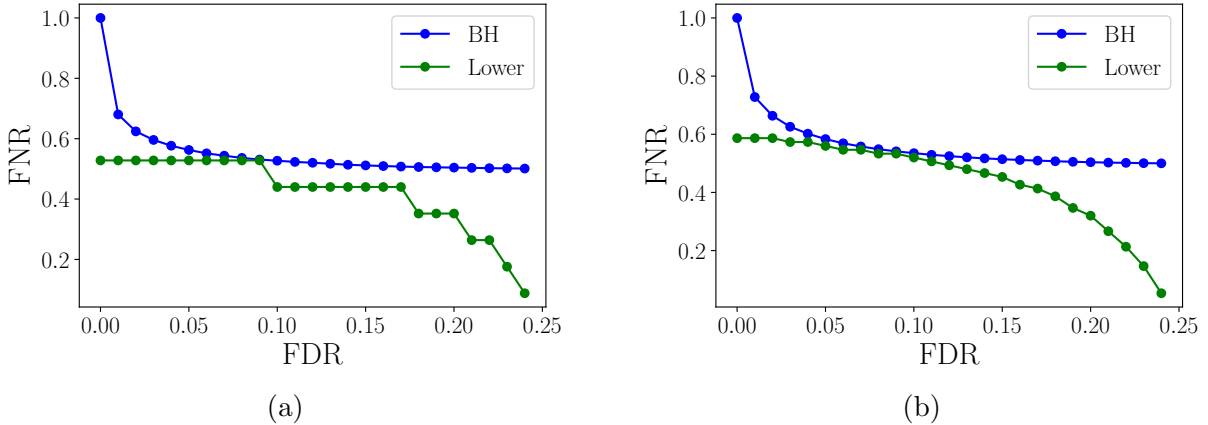

  \begin{center}
  \begin{tabular}{ccc}
\widgraph{0.5\linewidth}{{lower-bounds-curve-m15-r0.8}.pdf} & &
\widgraph{0.5\linewidth}{{lower-bounds-curve-m99-r0.6}.pdf} \\
(a) && (b)
  \end{tabular}
  \caption{The lower bound predicted by our theory plotted against the
    actual FDR-FNR tradeoff achieved by the Benjamini-Hochberg
    algorithm for two problems. FDR is on the horizontal axis, while
    FNR is on the vertical axis.  Both models use Gaussian test
    statistics with additive shifts $\mu = \sqrt{2r \log{n}}$ and $n =
    10000$. (a) Plots for $m = 15$ signals and $r = 0.8$.  (b) Plots
    for $m = 100$ and $r = 0.6$. Further details on these simulation
    experiments are given in
    Appendix~\ref{SecSimulation}. \label{FigPredictedCurves}}
  \end{center}
\end{figure}

In order to illustrate applications of this meta-theorem in practice,
we apply it to several specific models. These include the
previously-studied independent generalized Gaussians
model~\citep{Ari16DistFreeFDR,Rab20Tradeoff}, as well as Gaussian
models with dependence, Gaussian models with scale-transformation, and
an exponentiation model for $p$-values. Furthermore, we showcase our
approach's capacity to produce numerical lower bounds for concrete
models in Figure~\ref{FigPredictedCurves}.

\paragraph{Related work.}

Our work is most closely related to our own past work on the
non-asymptotic tradeoffs between FDR and
FNR~\citep{Rab20Tradeoff}. These previous results apply to a very
specific location shift model, in which the test statistics are
assumed to have tails on the order of $\exp\big(-|x|^{\gamma}\big)$
for some $\gamma \geq 1$---in other words, generalized Gaussian-like
tails. The tradeoffs derived in this past work, as well as previous
asymptotic lower bounds~\citep{Ari16DistFreeFDR}, apply only to the
independent case. There is another line of related work~\cite{Ji12UPS,
  Ji14Optimal}, applicable to models with exactly Gaussian noise,
that provides lower bounds that continue to hold under dependence, at
the cost of replacing FDR with the so-called \emph{modified} FDR
(mFDR), in which the expectation is moved inside the ratio in the
definition of FDR (cf.~\eqref{EqFDRFNR}). The mFDR and FDR measures
ought to behave similarly for large numbers of tests, but they are
distinct metrics, and the analysis strategies that work for mFDR
generally do not apply to FDR itself.

Our work is inspired in part by the work of Jin and Donoho on Tukey's
higher criticism (e.g.,~\citep{Don04HC,Don15HC,Jin14RareWeak}). Their
results apply to the Gaussian sequence model with location shifts, in
which the number of signals is a polynomially small fraction
$n^{-\betapar}$ of the total number of tests and in which signals are
weak (achieved by scaling the shift as $\mu = \sqrt{2r \log{n}}$ for
$0 < r < 1$). Their work establishes the regime of $\betapar$ and $r$
in which asymptotic consistency is possible, under the standard Type-I
and Type-II error measures for testing the global null.

\section{Background and problem formulation}

In this section, we provide necessary background and a precise formulation of
the problem under study.


\subsection{Multiple testing and false discovery rate}

Suppose that we observe a real-valued sequence $\Xseq \defn \{X_1,
\ldots, X_\numobs \}$ of $\numobs$ independent random variables.  
We introduce the
sequence of binary labels $\{H_1, \ldots, H_\numobs \}$ to encode
whether or not the null hypothesis holds for each observation; the
setting $H_i = 0$ indicates that the null hypothesis holds.  We define
\begin{align}
  \nulls & \defn \{ i \in [\numobs] \mid H_i = 0 \}, \quad \mbox{and}
  \quad \sigs \defn \{ i \in [\numobs] \mid H_i = 1 \}, \;
\end{align}
corresponding to the set of \emph{nulls} and \emph{signals},
respectively.  Our task is to identify a subset of indices that
contains as many signals as possible, while not containing too many
nulls.

More formally, a testing rule $\discoveries: \R^\numobs \rightarrow
2^{[\numobs]}$ is a measurable mapping of the observation sequence
$\Xseq$ to a set $\discoveries(\Xseq) \subseteq [\numobs]$ of
\emph{discoveries}, where the subset $\discoveries(\Xseq)$ contains
those indices for which the procedure rejects the null hypothesis.
There is no single unique measure of performance for a testing rule
for the localization problem.  In this paper, we study the notion of
the \emph{false discovery rate} (FDR), paired with the \emph{false
  non-discovery rate} (FNR).  These can be viewed as generalizations
of the type-I and type-II errors for single hypothesis testing.

We begin by defining the false discovery proportion (FDP), and false
non-discovery proportion (FNP), respectively, as
\begin{align}
\FDP_\numobs(\discoveries) \defn \frac{\card{\discoveries(\Xseq) \cap
    \nulls}}{ \card{\discoveries(\Xseq)} \vee 1}, \quad \mbox{and}
\quad \FNP_\numobs(\discoveries) \defn
\frac{\card{\discoveries(\Xseq)^{c} \cap \sigs}}{\card{\sigs}}.
\end{align}
Since the output $\discoveries(\Xseq)$ of the testing procedure is
random, both of these quantities are random variables.  The FDR and
FNR are given by taking the expectations of these random
quantities---that is
\begin{align}
  \label{EqFDRFNR}
\FDR_\numobs(\discoveries) & \defn \E\Big[\frac{ \card{
      \discoveries(\Xseq) \cap\nulls } }{\card {\discoveries(\Xseq)}
    \vee 1}\Big], \quad \mbox{and} \quad \FNR_\numobs(\discoveries)
\defn \E\Big[\frac{ \card{ \discoveries(\Xseq)^{c} \cap \sigs }
  }{\card{\sigs}} \Big],
\end{align}
where the expectation is taken over the random samples $\Xseq$.  Here
our definition of FNR follows\footnote{There is an alternative
  definition of $\FNR_{alt}$, in which the denominator is set to the
  number of non-rejections. In general, however, the number of
  non-rejections will be close to $n$ for any procedure with low
  $\FDR$ and thus in the sparse regime, the $\FNR_{alt}$ would
  trivially go to zero for any procedure that controls $\FDR$ at any
  level strictly below $1$.  The definition adopted here is therefore
  better suited to studying transitions in difficulty in the multiple
  testing problem.}  that of the paper~\cite{Ari16DistFreeFDR}.


\subsection{Model structure}
\label{SecModel}

In this paper, we consider a flexible class of models that includes as
special cases the location and scale families that have been studied
in past work.  For a testing problem with $n$ hypotheses, we assume
that the vector $X = (X_{1}, \dots, X_{n})$ of test statistics is
generated from some underlying random vector $W = (W_{1}, \dots,
W_{n})$ in the following way. The vector $W$ may be drawn from an
arbitrary (not necessarily product) distribution $\PW$, while the test
statistics are related via
\begin{align}
\label{EqGenModel}
X_{i} = \begin{cases} W_{i} &~\text{if}~i \in \nulls\\ \transf(W_{i})
  &~\text{otherwise.}
\end{cases}
\end{align}
The main restriction on this model is that $\transf \colon
\R\rightarrow\R$ must be a non-decreasing function such that
$\transf(w) \geq w$ for all $w$ in the support of a marginal of
$\PW$. For convenience, we shall also assume that the marginals of
$\PW$ are atom-free.  In some instances, we consider the restricted
case where all the $W_{i}$ are iid---that is, $\PW = \P_{0}^{\otimes
  n}$. We refer to this as the iid model.  Conceptually, a model is
therefore a tuple $(\PW,~\transf)$ satisfying these constraints, and
where appropriate we shall denote models by such tuples and sometimes
name them. \\

\noindent Prototypical examples of this general set-up include the
following:
\begin{description}
\item[Location model:] The variables $W_{i}$ are drawn from a
  generalized Gaussian distribution with density proportional to
  $\exp(-|w|^\gamma)$ for some $\gamma \in [1, 2]$, and the
  transformation function takes the form \mbox{$\transf(w) = \mu + w$}
  for some $\mu > 0$.
\item[Scale models:] The variables $W_{i}$ are the absolute values of
  standard normal variates, and the transformation function takes the
  form $\transf(w) = \sigma w$ for some $\sigma \geq 1$.
\item[Lehmann alternative model:] The variables $W_{i}$ are uniform on
  the unit interval, and the transformation function takes the form
  $\transf(w) = 1 - \big(1 - w\big)^{1/\lehmann}$ for some parameter
  $\lehmann \in (0,1)$.  This set-up models the situation in which the
  $W_{i}$ represent $p$-values and the signals have $p$-values that
  are stochastically closer to zero than those of the nulls. Since we
  have chosen to model the transformation as \emph{non-decreasing}, we
  represent the unit interval backwards, which leads to the form
  written rather than $w^{1/\lehmann}$.
\end{description}
We note that all three of these examples have been studied in past
work
(e.g.,~\citep{Ari16DistFreeFDR,Efr12LSI,Jin14RareWeak,Neu12FDRClass,Rab20Tradeoff}).

\subsection{\Topk procedures}

Many popular procedures, including the Benjamini-Hochberg (BH) and
several variants
thereof~\cite{Bar15Knockoffs,Ben95FDR,RamBarWaiJor19}, are based on
thresholding the order statistics.  Recall that the order statistics
of a sequence $\{X_1, \ldots, X_n\}$ are defined as
\begin{align}
\min_{i = 1, \ldots, n} X_i = \,: X_{(1)} \leq X_{(2)} \leq X_{(3)}
\leq \cdots \leq X_{(n)} \defn \max_{i = 1, \ldots, n} X_i.
\end{align}
A \emph{top-K procedure} is a method that rejects the hypotheses
corresponding to the top $K$ order statistics, where $K = K_n(\Xseq)$
is a non-negative integer that can depend on the observed statistics.
The testing rule $\discoveries \colon \R^{\numobs} \rightarrow
2^{[\numobs]}$ defined by any top-K procedure has the form
\begin{align}
\discoveries\big(\Xseq\big) = \bigg\lbrace i\in [\numobs] \colon X_{i}
\geq X_{(K_{n}(\Xseq))}\bigg\rbrace,
\end{align}
where $K_{n} \colon \R^{\numobs} \rightarrow \N$ is some (possibly
randomized) mapping.  Alternatively, such procedures can be described
in terms of choosing a threshold $\tau = \tau_n(\Xseq)$, and rejecting
all nulls $i$ for which $X_i \geq \tau$.


\section{Main results}
\label{SecMain}

We now turn to the statements of our main results.  We begin in
Section~\ref{SecDefProxy} by defining the deterministic proxies that
play a central role in our analysis; see Section~\ref{SecIntuition}
for the intuition that underlies these definitions.  In
Section~\ref{SecGeneral}, we state a general lower bound
(Theorem~\ref{ThmProxies}) on the pairs of FDR and FNR that are
achievable.  In the remaining sections, we illustrate the consequences
of this general bound for various specific models.


\subsection{A general bound based on deterministic proxies for FDR and FNR}
We say that a FDR-FNR pair $(\fdr, \fnr) \in [0,1] \times [0,1]$
is \emph{achievable} if there exists a \topk procedure $K$ such that
\begin{align}
  \label{EqnQcontrol}
  \FDR \big(K\big) \leq \fdr \quad \mbox{and} \quad \FNR\big(K\big)
  \leq \fnr.
\end{align}
Any \topk procedure satisfying condition~\eqref{EqnQcontrol} is said
to be \emph{$(\fdr, ~\fnr)$-controlled.}  Our high-level goal is to
provide bounds on the region of achievable $(\fdr, \fnr)$ pairs.


\subsubsection{Defining the deterministic proxies}
\label{SecDefProxy}

In order to characterize the space of achievable $(\fdr, \fnr)$ pairs,
we construct two sets of deterministic proxies. One set of proxies is
useful in the regime $\fdr \geq \fnr$, while the other is useful in
the opposite setting---namely, $\fnr \geq \fdr$.  The proxies used in
the former regime are denoted by $\FDPm$ and $\FNPm$, whereas those in
the latter regime are denoted by $\FNPp$ and $\FDPp$.  The reasoning
underlying our choice of notation should be clear once we detail their
construction below.  The proxies depend on the given pair $(\fdr$,
$\fnr)$, the model $\model = \big(\PW,\transf\big)$ under
consideration, and a parameter $\epsilon \in (0,1)$ that controls the
strength of the bounds. We make these dependencies explicit when
needed, suppressing them otherwise.

Our first step is to define deterministic approximations for the number
of total discoveries.  Letting $m = \card{\sigs}$ denote the number of
signals and given any $\epsilon > \max \{\fdr, \fnr \}$, we define
\begin{align}
  \label{EqnDiscProxies}
  \kstarm(\fnr,\epsilon) \defn \left ( 1 - \frac{\fnr}{\epsilon}
  \right ) m ~~ \text{and} ~~ \kstarp(\fdr,\epsilon) \defn \left ( 1 -
  \frac{\fdr}{\epsilon} \right)^{-1} m .
\end{align}
Roughly speaking, the integer $\kstarm$ functions as a lower
approximation for the number of \emph{total discoveries}, whereas the
quantity $\kstarp$ provides an upper approximation for the same
quantity.  Note that these lower and upper bounds converge as $\fdr,
\fnr \rightarrow 0$; in the limit $\fdr = \fnr = 0$, we have
$\kstarm(0, \epsilon) = \kstarp(0, \epsilon) = m$, since in this case,
the total number of discoveries must be equal to the number of signals
$m$.

For each of these approximations of the number of discoveries, we construct
a corresponding false discovery proxy.  Recalling the 
random vector \mbox{$W = (W_{\nulls}, W_{\sigs})$} that underlies our
generic model, these quantities involve the order statistics
\begin{align*}
W_{\nulls, (1)} \leq W_{\nulls, (2)} \leq \cdots \leq W_{\nulls,
  (|\nulls|)},
\end{align*}
with the order statistics for $W_{\sigs}$ defined analogously.  For
any $\epsilon > \max \{ \fdr, \fnr \}$, adopting the shorthand
$\kstarm = \kstarm(\fnr, \epsilon)$ and $\kstarp = \kstarp(\fdr,
\epsilon)$, we define proxies as follows.
\paragraph{False discovery proxies:}
\begin{subequations}
      \label{EqnFalseProxy}
  \begin{align}
\lstarm(\fnr,\epsilon,\model) & = \arg \max_{\ell \in [1, \kstarm] }
\left \{ \P \Big[ W_{\nulls,(\ell)} > \transf\big(W_{\sigs,(\kstarm -
    \ell + 1)}\big) \Big] \geq 1 - \epsilon \right \}, \qquad
\mbox{and} \\
\lstarp(\fdr,\epsilon,\model) & = \arg \max_{\ell \in [ \kstarp -m,
    \kstarp]} \left \{ \P \Big[ W_{\nulls,(\ell)} >
  \transf\big(W_{\sigs,(\kstarp - \ell + 1)}\big) \Big] \geq 1 -
\epsilon \right \}.
\end{align}
\end{subequations}
Roughly, the quantities $\lstarm$ and $\lstarp$ represent,
respectively, lower and upper approximations to the number of false
discoveries.

Finally, by taking appropriate ratios, we define:
\paragraph{Proxies to FDR and FNR:}
\begin{subequations}
\begin{align}
\FDP_{-}^{\ast}(\fnr,\epsilon,\model) =
\frac{\lstarm(\fnr,\epsilon,\model)}{m} ~~ \text{and} ~~
\FNP_{-}^{\ast}(\fnr,\epsilon,\model) = \frac{m -
  \kstarm(\fnr,\epsilon) + \lstarm(\fnr,\epsilon,\model)}{m} ,
\\ \FDP_{+}^{\ast}(\fdr, \epsilon,\model) = \frac{\lstarp(\fdr,
  \epsilon,\model)}{m} ~~ \text{and} ~~ \FNP_{+}^{\ast}(\fdr,
\epsilon,\model) = \frac{m - \kstarp(\fdr, \epsilon) + \lstarp(\fdr,
  \epsilon,\model)}{m} .
\end{align}
\end{subequations}
To be clear, in defining $\FDP_{+}^{\ast}$ (respectively
$\FDP_{-}^\ast$), it might be more natural to use $\kstarp$
(respectively $\kstarm$) in the denominator, but as noted above, when
$(\fdr, \fnr)$ are small, both of these quantities are close to $m$.


\subsubsection{The underlying intuition}
\label{SecIntuition}

Let us now describe the intuition that underlies the
definitions~\eqref{EqnFalseProxy}, ignoring the difference between the
$+$ and $-$ versions so as to simplify matters.  First, suppose that
we accept that $k$ is a good approximation to the total number of
discoveries, and that $\ell$ is a good approximation to the number of
false discoveries.  In this case, $\ell/k$ is a good approximation to
the FDR, and since $k - \ell$ of the discoveries must be false, then
$\frac{m - (k - \ell)}{m}$ should be a good approximation to the FNR.
As we have argued above, when $\fnr$ and $\fdr$ are small, then $k$ is
actually relatively close to $m$, so that $\ell/m$ should also be a
good approximation to the FDR.

It remains to justify why $\ell$, as defined in
equation~\eqref{EqnFalseProxy}, is a reasonable proxy to the number of
false discoveries.  Consider a procedure that rejects exactly $k$
hypotheses, of which $\ell$ are nulls.  It must then be case that
$\ell^{\text{th}}$ largest null value exceeds the value of the $\big(k
- \ell + 1\big)^{\text{th}}$ largest signal value, or else only $\ell
- 1$ nulls would be in the top $k$ test statistics. Using the
definition of the model, we can re-express this relation in symbols:
\begin{align}
\label{EqOrdStatsRelation}
W_{\nulls,(\ell)} > \transf \left(W_{\sigs,(k - \ell + 1)} \right).
\end{align}
The definitions~\eqref{EqnFalseProxy} are motivated by this assertion.


\subsubsection{A general lower bound}
\label{SecGeneral}

Our main result is that our choice of proxies yield constant-factor
lower bounds on the attainable FDR and FNR of any \topk procedure.
\begin{thm}
\label{ThmProxies}
Given a model $\model$, consider any $(\fdr,~\fnr)$-controlled \topk
procedure such that \mbox{$2 \max \{ \fdr, \fnr \} < \frac{1}{3}$.}
Then for any scalar $\epsilon \in \big( 2 \max \{\fdr, \fnr \},
\frac{1}{3} \big)$, there exists a constant $\cproxy(\epsilon) \geq 1$
such that
\begin{subequations}
  \begin{align}
    \label{EqnLineOne}
\fdr \geq \cproxy^{-1} \FDPm(\fnr,\epsilon) ~~ & \text{and}
~~\max\big\lbrace \fdr ,~\fnr\big\rbrace \geq
\cproxy^{-1}\FNPm(\fnr,\epsilon), \qquad \mbox{as well as } \\
\label{EqnLineTwo}
\fnr \geq \cproxy^{-1}\FNPp(\fdr,\epsilon) ~~ & \text{and}~~
\max\big\lbrace \fdr, ~\fnr\big\rbrace \geq \cproxy^{-1}\FDPp(\fdr,
\epsilon).
\end{align}
\end{subequations}
\end{thm}
The slightly unorthodox form of~\eqref{EqnLineOne} and~\eqref{EqnLineTwo} calls for some discussion. The presence of the maximum reflects the fact that generally only one set of proxies
will be suitable for lower bounding $\fdr$ and $\fnr$ simultaneously. If $\fdr > \fnr$, the bound in~\eqref{EqnLineTwo} is the meaningful one, while~\eqref{EqnLineOne} gives the desired bound when $\fnr > \fdr$. When
$\fdr = \fnr$, either equation will do. 

The two regimes arise because $\lstarp$ and $\kstarp$ yield a good approximation of the false discovery number and total number of discoveries only when $\fdr > \fnr$, while $\lstarm$ and $\kstarm$ provide a good approximation only when $\fdr < \fnr$. Intuitively, the dichotomy arises because $\kstarp$ may be \emph{larger} than the actual number of discoveries by an amount as large as order of $\fdr + \fnr$, so
that $\FDPp$ can only be upper bounded by a quantity of this order, or, equivalently (disregarding constants), a quantity on the order of $\max\left\{ \fdr, \fnr\right\}$. A similar but inverted phenomenon occurs for the $-$ proxies. 


\subsection{Application 1: Independent Gaussians model}

In this section, we investigate models in which the vector $W$ has iid
Gaussian entries, and the signal structure is specified by either a
location shift or a scale factor.


\subsubsection{Gaussian location model}

We begin by analyzing the Gaussian location model, in which the
function $\transf$ takes the form
\begin{align}
  \transf(w) & = w + \mu \qquad \mbox{for some $\mu > 0$. }
\end{align}
By applying Theorem~\ref{ThmProxies} to this particular model, we
obtain the following:
\begin{cor}
  \label{ThmIIDNormal}
  Consider the iid Gaussian location model with \mbox{$m = n^{1 -
      \betapar}$} signals and \mbox{$\mu = \sqrt{2r \log{n}}$} with
  parameters $(\betapar, r)$ satisfying the inequality $0 < \betapar <
  r < 1$. Suppose that there exists a constant $c > 0$ such that for
  all $n \geq 1$, there is an $(\fdr_n,~ \fnr_{n})$-controlled top-$K$
  procedure with $\fdr_{n} = cn^{-\kappa_\fdr}$ and $\fnr_{n} =
  cn^{-\kappa_\fnr}$.  Then we must have
  \begin{align}
    \label{EqnIIDNormalBound}
    \sqrt{\betapar + \kappa_\fdr} + \sqrt{\kappa_\fnr} \leq \sqrt{r},
    \quad \mbox{and} \quad \min\big \lbrace \kappa_\fdr, ~\kappa_\fnr
    \big\rbrace \leq \underbrace{\frac{\big(r - \betapar
      \big)^{2}}{4r}}_{= : \; \kappa^\ast}.
\end{align}
\end{cor}

The result obtained in~\eqref{EqnIIDNormalBound} is essentially the
lower bound of Rabinovich et al.~\cite{Rab20Tradeoff}, derived by
other means in that paper and applicable to an exactly Gaussian rather
than a Gaussian-like model. Thus, Corollary~\ref{ThmIIDNormal} can be
seen as an extension of those earlier results, illustrating how the
methods developed in this work can expand the scope of multiple
testing lower bounds. Moreover, since the lower bound of Rabinovich et
al.~\cite{Rab20Tradeoff} is known to be sharp, it seems likely the
bound proven here is likewise sharp for the Gaussian location model.

Figure~\ref{FigPredictedCurves} shows the predictions of
Theorem~\ref{ThmProxies} in a Gaussian model, in particular giving
plots of the lower bound predicted by our theory plotted against the
actual FDR-FNR tradeoff achieved by the Benjamini-Hochberg algorithm
for two problems. In each plot, the FDR is on the horizontal axis,
while FNR is on the vertical axis.  See the figure caption for further
details.


\subsubsection{Gaussian scale model}

We now turn attention to the Gaussian scale model.  It is specified by
the transformation
\begin{align*}
\transf(w) & = \sigma w, \qquad \mbox{for some $\sigma \geq 1$.}
\end{align*}
By applying Theorem~\ref{ThmProxies} to this model, we obtain a rather
different lower bound on pairs $(\fdr, \fnr)$.  At a high level, the
main take-away is that the FDR and FNR can only decay as inverse
polynomial functions of $n$ when the signal strength $\sigma$ is
extremely strong---in particular, the scalar $\sigma$ has to grow
polynomially in $n$.

\begin{cor}
\label{ThmIIDNormalScale}
Consider the iid Gaussian scale model with $m$ signals and signal
strength $\sigma \geq 1$ where $\betapar_{n} \defn \frac{\log
  \frac{n}{m}}{\log{n}}$ lies in the interval $[\rho,~1 - \rho]$ for
some $\rho \in (0, 0.5)$.  Suppose that there exists an $(\fdr,
~\fnr)$-controlled procedure such that $\max\big\lbrace \fdr,
~\fnr\big\rbrace \leq \frac{1}{3\cproxy}$.  Then there exists some
$\eta_{n} \in (0,1)$ such that
\begin{align}
\sigma & \geq \frac{1}{\sqrt{2\pi}\cproxy} \cdot \bigg(1 -
\eta_{n}\bigg)\bigg(\frac{1}{m} + \fnr\bigg)^{-1}\sqrt{2\betapar_{n}
  \log{n} + 2\log\bigg(\fdr + \frac{1}{m}\bigg)^{-1}}.
\end{align}
\end{cor}

As in the statement of Corollary~\ref{ThmIIDNormal} for the location
model, by assuming certain scalings of the number of signals, FDR and
FNR, we can give an asymptotic statement.  In particular, suppose that
the number of signals scales as $m \propto n^{1 - \betapar}$ for a
fixed $\betapar$, whereas the FDR and FNR scale as $\fdr_{n} \propto
n^{-\kappa_\fdr}$ and $\fnr_{n} \propto n^{-\kappa_\fnr}$ for some
scalars $\kappa_\fdr, \kappa_\fnr$ such that \mbox{$\max
  \lbrace\kappa_\fdr, ~\kappa_\fnr \rbrace \leq 1 - \betapar$.}  Then
there is a universal constant $c > 0$ such that
\begin{align}
\sigma \geq c \; n^{\kappa_\fnr} \sqrt{2\big(\betapar + \kappa_\fdr
  \big)\log{n}}.
\end{align}
Consequently, we see that whenever $\kappa_\fnr > 0$, the signal
strength $\sigma$ must grow polynomially in $n$.  This is quite a
dramatic contrast from the location model, where the analogous
quantity $\mu$ need only grow proportionally to $\sqrt{\log n}$.


\subsection{Application 2: Gaussian location models with dependence}

Given the presence of dependence in many target applications of
multiple testing (e.g.,~\citep{Bar16PFilter}), it makes sense to ask how dependence
changes the performance of multiple testing procedures. In this
section, we provide answer for two models of Gaussian dependence
which lie at opposite extremes of dependency.  In both
cases, we consider only location shifts.


\subsubsection{Spiked dependence model}

We begin by considering a doubly-spiked covariance model, with one
spike within the nulls and a separate spike within the signals.  This
specification corresponds to coupling all the nulls (and, separately,
all the signals) through a single random variable that captures all
shared randomness.

In the spiked dependence model, the model distribution $\PW$ is a
multivariate Gaussian $\Norm\big(0,~\Sigma\big)$, with covariance
matrix in the block-partitioned form
\begin{align}
\label{EqSpikedModel}
\Sigma_{ij} = \begin{cases} 1 &~\text{if}~i = j, \\
  \rhonull &~\text{if}~i \neq j ~~ \text{and} ~~ i, j \in \nulls, \\
  \rhosig &~\text{if}~i \neq j~~ \text{and} ~~ i, j \notin \nulls, \\
  \pm \rhocross &~\text{if}~i \neq j ~~ \text{and}~ i \in \nulls,~j
  \notin \nulls ~~ \text{or vice versa},
\end{cases}
\end{align}
for parameters $0 \leq \rhonull, \rhosig, \rhocross < 1$.

\begin{cor}
\label{ThmSpikedNormal}
Consider the spiked dependence model with $m = n^{1 - \betapar}$ and
$\mu = \sqrt{2r \log{n}}$ for some pair $(\betapar, r)$ satisfying the
inequalities $0 < \betapar < r < 1$.  Suppose that for each $n$, there
exists an $(\fdr_n,~ \fnr_{n})$-controlled top-$k$ procedure with
$\fdr_{n} = c n^{-\kappa_\fdr}$ and $\fnr_{n} = c n^{-\kappa_\fnr}$.
Then we must have
\begin{align}
  \label{EqnSpikedNormalBound}
\sqrt{1 - \rhonull} \sqrt{\betapar + \kappa_\fdr} + \sqrt{1 -
  \rhosig}\sqrt{\kappa_\fnr} \leq \sqrt{r}.
\end{align}
\end{cor}

Note that the bound~\eqref{EqnSpikedNormalBound} is a generalization
of the bound~\eqref{EqnIIDNormalBound} for iid Gaussians, to which it
reduces when $\rhonull = \rhosig = \rhocross = 0$.  Relative to this
iid case, the bound~\eqref{EqnSpikedNormalBound} allows for larger
values of the pair $(\kappa_\fdr, \kappa_\fnr)$---which translates
into faster decay rates of FDR and FNR---when either $\rhosig$ or
$\rhonull$ is non-zero.  While this might be counterintuitive at first
sight, note that our spiked dependence makes all nulls more similar to
each other when $\rhonull > 0$, and all signals more similar to each
other when $\rhosig > 0$.  This similarity in either the nulls or
signals means that it becomes easier to control the FDR and FNR.  What
may still be surprising is that $\rhocross$ does not play any role in
the rates.


\subsubsection{Grouped dependence model}

We now turn to the opposite extreme of dependency.  In the
grouped dependence model, we match each signal with a different set of
$A$ nulls that are strongly coupled to that signal, but independent of
all other signals and all nulls in different groups. More formally,
first fix a value $1 \leq A \leq \min\big\lbrace
m,~\frac{n}{m}\big\rbrace$. We then write $\sigs = \lbrace i_{1} <
\cdots < i_{m}\big\rbrace$ and for each $1 \leq g \leq m$, we define a
set of $A$ nulls $\nulls^{(g)}$ corresponding to the $g^{\text{th}}$
signal. Finally, we define the independent nulls as $\nulls^{(0)} =
\nulls\setminus \cup_{g = 1}^{m} \nulls^{(g)}$.

Rather than providing an explicit form of the covariance matrix for
this model, it is more informative to specify the underlying
generative model, given by
\begin{align}
\label{EqGroupedModel}
W_{\sigs} \sim \Norm\big(0,~I_{m}\big), \quad \mbox{with} \quad ~~
W_{i} ~|~W_{\sigs} \sim \begin{cases}
  \Norm\big(0,~1\big) &~\text{if}~ i \in \nulls^{(0)}, \\
  W_{i_g} &~\text{if}~i \in \nulls^{(g)}.
\end{cases}
\end{align}
By applying Theorem~\ref{ThmProxies} to this model, we obtain the
folllowing:

\begin{cor}
\label{ThmGroupedNormal}
Consider the grouped dependence model with $m = n^{1 - \betapar}$,
$\mu = \sqrt{2r \log{n}}$, and $A = \frac{n - m}{n^{1 - \etapar}}$
where the parameters $(\betapar, r, \etapar)$ satisfy the inequalities
\mbox{$0 \leq \etapar < \betapar < r < 1$.}  Suppose that there is a
constant $c > 0$ such that for each positive integer $n$, there is an
\mbox{$(\fdr_{n},~\fnr _{n})$-controlled} procedure with $\fdr_{n} = c
n^{-\kappa_\fdr}$ and $\fnr_{n} = c n^{-\kappa_\fnr}$. Then
\begin{align}
  \label{EqnGroupedNormal}
\sqrt{\betapar + \kappa_\fdr} + \sqrt{\kappa_\fnr} \leq \sqrt{r}.
\end{align}
\end{cor}

Note the surprising fact that the
bound in~\eqref{EqnGroupedNormal} is identical to our earlier
bound~\eqref{EqnIIDNormalBound} from the iid case. If the lower bound of~\eqref{EqnGroupedNormal}
is sharp, this coincidence reflects a deep fact about the difficulty of multiple testing 
in the grouped Gaussians model. We do not at this point know, however, whether
Corollary~\ref{ThmGroupedNormal} is sharp, and the sharpness of this result (like the others established in this paper)
remains an important question for future work. 



\subsection{Application 3: Lehmann alternatives}

The Lehmann alternative model has often been used in theoretical analyses of multiple-testing procedures. 
In this model, the statistics are now $p$-values; nulls are assumed to come from a
uniform distribution, while alternatives follow a CDF given by
\begin{align}
\label{eq:lehmann-alt}
F(p) & = p^{\lehmann} \qquad \mbox{for some $\lehmann \in (0,1)$.}
\end{align}

In order to formulate this problem within our framework, let $W_{i}$
be iid uniform random variables on the unit interval $[0,~1]$, and
define the transformation
\begin{align*}
\transf(w) = 1 - \big(1 - w\big)^{1/\lehmann} .
\end{align*}
Note that here $1 - w$ plays the role of the $p$-value, so that the
$w$ values for signals are more clustered around $1$ than is the case for the nulls.

\begin{cor}
\label{ThmLehmannAlt}
Under the Lehmann alternative model with parameter $\lehmann \in (0,1)$,
fix some triple $(\fdr, \fnr, \epsilon)$ such that $\fdr \leq
\frac{\epsilon}{3}$ and $\etapar \defn \frac{3 \fnr}{\epsilon} +
\frac{1}{m} + \sqrt{\frac{3 \cproxy \fnr}{m \epsilon}} < 1$, where $\cproxy$ is the constant from Theorem~\ref{ThmProxies}.  
Further, let $\pi_{1} = \frac{m}{n}$. Then for any $(\fdr, \fnr)$-controlled procedure, we must have
\begin{align}
\frac{1}{\lehmann} \geq \frac{1 - \etapar}{\etapar} \cdot \log
\Bigg(\frac{\epsilon}{3 \pi_{1} \fdr} \bigg[1 + 4
  \log{\frac{3}{\epsilon} }\bigg]^{-1}\Bigg).
\end{align}
\end{cor}


The bound of Corollary~\ref{ThmLehmannAlt} requires some interpretation. Intuitively, $t$ is on the order of $\fnr$, while the argument of the logarithm is on the order of $\frac{1}{\pi_{1}\fdr}$, so the high-level takeaway is that
\begin{equation}\label{EqLehmannAltSimple}
\frac{1}{\lehmann} \gtrsim \frac{1}{\fnr} \log \frac{1}{\pi_{1}\fdr}.
\end{equation}
Since the signal---the difference between nulls and alternatives---becomes greater as $\lehmann$ becomes smaller, $\frac{1}{\lehmann}$ is a measure of signal strength, 
and thus~\eqref{EqLehmannAltSimple} is similar to our previous bounds in that it lower bounds the required signal strength in terms of the problem parameters. In this case, the dependence on
the FNR $\fnr$ is inverse polynomial, while the dependence on both the FDR $\fdr$ and the sparsity $\pi_{1}$ is logarithmic.


\section{Proofs}

We now turn to the proofs of our results.
\subsection{Technical tools}

Before giving proofs of our main results, we develop two technical
tools that we apply repeatedly in our arguments.


\subsubsection{Derandomization under concentration}

Our proxies depend on the model only through the false discovery
number proxies $\lstarm$ and $\lstarp$. Unfortunately, the dependence
is of a rather complicated form, since the
definitions~\eqref{EqnFalseProxy} involve the probabilities of events
defined in terms of the order statistics of nulls and signals. In
order to make progress, our first step is to simplify these
definitions so as to obtain modified versions that are more tractable.
In this section, we show that, provided the model's order statistics
admit a suitable concentration bound, we can reduce the probabilistic
comparison to a deterministic comparison of expected order statistics.
In particular, we make use of the following family of assumptions,
which are parameterized by $T \in \lbrace \nulls, \sigs\rbrace$, and
an integer $k$.

\paragraph{Concentration assumption $(T, k)$:}
There exists a function $\ordD_{T,k} \colon (0,~1) \rightarrow
[0,~\infty)$ such that
\begin{align}
\label{EqnConcentration}
\P \Biggr [ \big|X_{T, (k)} - \E\big[X_{T,(k)}\big]\big| \geq
  \ordD_{T,k}\big(\epsilon\big) \Biggr] \leq \epsilon.
\end{align}
Depending on the nature of the function $\ordD_{T, k}$,
condition~\eqref{EqnConcentration} might be a more or less stringent
(and a more or less useful) assumption.  In general, when we apply
this bound, we shall be able to prove it holds for a reasonable choice
of $\ordD_{T,k}$. \\

\noindent Our analysis invokes two particular cases of the
concentration assumption:
\begin{description}
\item[Case I:] The concentration assumption~\eqref{EqnConcentration}
  holds for $(\nulls, \lstar)$ and $(\sigs, \kstar - \lstar + 1)$.
\item[Case II:] The concentration assumption~\eqref{EqnConcentration}
  holds for $(\nulls, \lstar+1)$ and $(\sigs, \kstar - \lstar)$.
\end{description}
Here as always, the integer $\kstar$ is one of $\kstarm$ and
$\kstarp$, and the integer $\lstar$ is fixed correspondingly.  The following
lemma allows us to reduce from probabilities of events to differences
in expected order statistics.

\begin{lem}
\label{LemDerandom}
Suppose that $\lstar \in \{ \lstarm\big(\fnr, \epsilon\big),
\lstarp\big(\fdr, \epsilon\big) \}$ and define $\kstar$
accordingly. Then under Cases I and II, we have the following bounds:
\begin{subequations}
\begin{align}
\label{EqDerandom}
{\bf{(I):}} \quad \E\bigg[\transf\big(W_{\sigs,(\kstar -
    \lstar)}\big)\bigg] + \ordD_{\sigs,(\kstar -
  \lstar)}\bigg(\frac{\epsilon}{3}\bigg) > \E\bigg[W_{\nulls,(\lstar +
    1)}\bigg] - \ordD_{\nulls,(\lstar +
  1)}\bigg(\frac{\epsilon}{3}\bigg). \\
\label{EqDerandomConv}
{\bf{(II):}} \quad \E\bigg[\transf\big(W_{\sigs,(\kstar - \lstar +
    1)}\big)\bigg] - \ordD_{\sigs,(\kstar - \lstar +
  1)}\bigg(\frac{\epsilon}{3}\bigg) <
\E\bigg[W_{\nulls,(\lstar)}\bigg] +
\ordD_{\nulls,(\lstar)}\bigg(\frac{\epsilon}{3}\bigg).
\end{align}
\end{subequations}
\end{lem}
\noindent See Appendix~\ref{AppLemDerandom} for the proof of this
claim. \\

\noindent {\bf{Remarks:}} The main value of Lemma~\ref{LemDerandom}
lies in inequality~\eqref{EqDerandom}. Indeed, this inequality places
a lower bound on an expected order statistic coming from a signal in terms of an
expected order statistic coming from a null (plus some deviations). Since signals
are shifted rightward relative to nulls, a lower bound of this kind
gives a lower bound on the signal strength in terms of
$\lstar$. Meanwhile, Theorem~\ref{ThmProxies} provides upper bounds on
$\lstar$ (and $m - \kstar + \lstar$) in terms of the realized FDR and
FNR. Together, these bounds yield a lower bound on signal strength in
terms of FDR and FNR, which can be interpreted as a lower bound on FDR
and FNR in terms of the signal strength.


\subsubsection{Transfering results between models}

It is useful to be able to transfer results from simple models to more
complex models that are in some sense ``close'' to them. In this
section, we specify a notion of closeness that makes sense for our
problem and prove a technical result that allows us to transfer lower
bounds between close models.

Our definition of closeness for models has some unusual features that
bear explanation. First, it only applies to models that share a single
transformation function $\transf$. This limitation is imposed for
convenience and is not fundamental. The definition is also asymmetric,
with some base model $\model$ given, and the closeness of another
model $\model'$ assessed relative to $\model$. The asymmetry arises
from the fact that we wish to define proximity of models based on a
single fixed distribution over the order statistics, rather than uniformly over all distrubtions over order statistics, and the single fixed reference
point we choose arises from the discovery and false discovery number
proxies $\kstar_{\pm}$ and $\lstar_{\pm}$, which depend on the model.

\begin{defenv}
\label{DefCloseness}
We say that two models $\model = \big(\PW,~\transf\big)$ and $\model'
= \big(\PW',~\transf\big)$ are $(\nullDelta, \sigDelta, \lstarp,
\delta)$-close if
\begin{align*}
\max\Bigg\lbrace \P\big(\big|W_{\nulls,(\lstar)} -
W_{\nulls,(\lstar)}'\big| \geq \nullDelta
\big),~\P\big(\big|f\big(W_{\sigs,(\kstar - \lstar)}\big) - \transf
\big( W_{\sigs,(\kstar - \lstar)}'\big) \big| \geq \sigDelta \big)
\Bigg\rbrace \leq \delta,
\end{align*}
where $\lstar = \lstarp(\fdr, \delta, \model)$ and $\kstar =
\kstarp(\fdr, \delta)$.  Similarly, we say that they are $(\nullDelta,
\sigDelta, \lstarm, \delta)$-close if the same condition holds with
$\lstar = \lstarm(\fnr, \delta, \model)$ and $\kstar = \kstarm(\fnr,
\delta)$.
\end{defenv}

Based on this definition, we can transfer lower bounds from the base
model in Definition~\ref{DefCloseness} to the other model using the
following technical lemma.

\begin{lem}
\label{LemTransfer}  
 For a given pair $(\fdr, \fnr)$ with $2 \max \{\fdr, \fnr \} <
 \frac{1}{3}$, consider some $\epsilon \in \big (2 \max \{\fdr, \fnr
 \}, \frac{1}{3} \big)$.
 \begin{enumerate}
\item[(a)] If $\big(\PW',~\transf\big)$ is $(\nullDelta,\sigDelta,
  \lstarp, \epsilon/3)$-close to $\big(\PW,~\transf\big)$, then
  \begin{subequations}
\begin{align}
 \lstarp \big(\fdr,~\epsilon,~\big(\PW',~f - \nullDelta -
 \sigDelta\big) \big) & \geq \lstarp(\fdr, \epsilon, (\PW, f)).
\end{align}
\item[(b)] If $\big(\PW',~\transf\big)$ is $(\nullDelta,\sigDelta,
  \lstarm, \epsilon/3)$-close to $\big(\PW,~\transf\big)$, then
\begin{align}
 \lstarm \big(\fnr,~\epsilon,~\big(\PW',~f - \nullDelta -
 \sigDelta\big) \big) & \geq \lstarm(\fnr, \epsilon, (\PW, f)).
\end{align}
  \end{subequations}
 \end{enumerate}
\end{lem}
\noindent See Appendix~\ref{AppLemTransfer} for the proof of this claim.

We use Lemma~\ref{LemTransfer} primarily to remove dependence. In that
context, a particularly useful specialization of it is the following
decoupling lemma, which allows us to remove dependence between nulls
and signals provided that we can verify the concentration
condition~\eqref{EqnConcentration}.

\begin{lem}
  \label{LemDecouple}
For a given model $\model = \big(\PW,~\transf\big)$, let $\model'$
denote the same model but with nulls and signals sampled independently
from their marginals under $\PW$. Suppose that $2 \max \{ \fdr, \fnr
\} < \frac{1}{3}$, and that $\model'$ satisfies Case I of the
concentration assumption~\eqref{EqnConcentration}.  Then for any
\mbox{$\epsilon \in \big( 2 \max \{\fdr, \fnr \}, \frac{1}{3} \big)$,}
with the integers $\kstar = \kstarp(\fdr, \epsilon)$, \mbox{${\lstar}'
  \defn \lstarp(\fdr, \model', \frac{\epsilon}{3})$,} the scalar
$\Delta = 2\big[\ordD_{\nulls,({\lstar}')}\big(\epsilon/6\big) +
  \ordD_{\sigs,(\kstar - {\lstar}' + 1)}\big(\epsilon/6\big)\big]$,
and the model \mbox{$\model'' = \big(\PW,~w \mapsto f(w) -
  2\Delta\big)$,} we have
\begin{subequations}
\begin{align}
\lstarp(\fdr, \model'', \epsilon) & \geq \lstarp(\fdr, \model',
\frac{\epsilon}{3}).
\end{align}
Similarly, with $\kstar = \kstarm(\fnr, \epsilon)$, \mbox{${\lstar}'
  \defn \lstarm(\fnr, \model', \frac{\epsilon}{3})$,} and $(\Delta,
\model'')$ redefined accordingly, we have
\begin{align}
\lstarm(\fnr, \model'', \epsilon) & \geq \lstarm(\fnr, \model',
\frac{\epsilon}{3} ).
\end{align}
\end{subequations}
\end{lem}
\noindent See Appendix~\ref{AppLemDecouple} for the proof of this claim.

\subsection{Proof of Theorem~\ref{ThmProxies}}

In this section, we prove Theorem~\ref{ThmProxies}. Two main ideas
underlie the proof. First, we show that any \topk procedure that is
$(\fdr, \fnr)$-controlled must have at least a constant probability of
making approximately the correct number of discoveries (meaning that
the number of discoveries is equal to the true number of signals).  In
order to formalize this idea, for a given \topk procedure $K$, define
the event
\begin{align}
\label{EqBandEvent}
\bandevent & \defn \eventwrap{K \in \bigg[\kstarm(\fnr),~\kstarp(\fdr)
    \bigg]} ,
\end{align}
where $\kstarm$ and $\kstarp$ are the discovery
proxies~\eqref{EqnDiscProxies}.  The width of this band is determined
by $\fdr$ and $\fnr$ and by the constant $\epsilon$, which will play
the role of a parameter in the analysis throughout our proofs.
\begin{lem}
\label{LemBand}
 For any $(\fdr,~\fnr)$-controlled \topk procedure, we have $\P
 [\bandevent ] \geq 1 - 2\epsilon$.
\end{lem}
\noindent We defer the proof of this lemma to Section~\ref{SecLemBand}. \\

The main second ingredient is a precise version of the argument that led
to the inequality~\eqref{EqOrdStatsRelation}. Essentially, we need to
know that the event defined by~\eqref{EqOrdStatsRelation} really is
the same as the event defined by the number of false discoveries in
the top $k$ being lower bounded by $\ell$.  We define
\begin{align}
\label{EqFalseDiscNotation}
\falsedisc{k} & = \bigg|\big\lbrace i \colon X_{i} \geq
X_{(k)}\big\rbrace \cap \nulls\bigg|,
\end{align}
corresponding to the number of false discoveries in the top $k$.
In terms of this notation, we have the following:
\begin{lem}
  \label{LemOrdStatsRelation}
We have
\begin{align}
\eventwrap{\falsedisc{k} \geq \ell} = \eventwrap{W_{\nulls,(\ell)} >
  \transf\big(W_{\sigs,(k - \ell + 1)}\big)} \qquad \mbox{for each $k
  = 1, 2, \ldots, n$.}
\end{align}
\end{lem}
\noindent See Section~\ref{SecLemOrdStatsRelation} for the proof of
this claim.\\

Equipped with these lemmas, we now turn to the proof of the theorem.
Define the events
\begin{align*}
\proxyeventm = \eventwrap{L(\kstarm(\fnr)) \geq \lstarm(\fnr)} ~~
\text{and} ~~ \proxyeventp = \eventwrap{L(\kstarp(\fdr)) \geq
  \lstarp(\fdr)}.
\end{align*}
By applying Lemma~\ref{LemOrdStatsRelation} twice, once with the
choice $k = \kstarm(\fnr)$ and then with the choice $k =
\kstarp(\fdr)$, and using the definitions of~$\lstarm$ and $\lstarp$
(see equation~\eqref{EqnFalseProxy}), we have
\begin{align}
\label{EqProxyEventProb}
\min \Big \{ \P [\proxyeventm] , ~\P [\proxyeventp] \Big \} & \geq 1 -
\epsilon.
\end{align}
Next, combining Lemma~\ref{LemBand} and the
bound~\eqref{EqProxyEventProb} yields
\begin{align}
\label{EqGoodEvents}
\P\bigg(\overbrace{\bandevent \cap \proxyeventm}^{\goodeventm}\bigg)
\geq 1 - 3\epsilon ~~ \text{and} ~~ \P\bigg(\overbrace{\bandevent \cap
  \proxyeventp}^{\goodeventp}\bigg) \geq 1 - 3\epsilon .
\end{align}

\paragraph{Argument for $-$ proxies:}  We now use the bound~\eqref{EqGoodEvents}
to proof the theorem's claims for the negative-subscript proxies.
Note that $L(K) = K \cdot \FDP(K)$, so that on conditioned on
$\goodeventm$, we have
\begin{align*}
& \lstarm \leq K \cdot \FDP(K) \leq \kstarp \cdot \FDP(K) \\ &
  ~~~~~~~\implies \FDP(K) \geq \frac{\lstarm}{\kstarp} = \left(1 -
  \frac{\fdr}{\epsilon} \right) \cdot \FDPm .
\end{align*}
We now take expectations to find that
\begin{align*}
\FDR(K) & \geq \P \big[ \goodeventm \big] \cdot \E\big[\FDP(K)~ \mid
  ~\goodeventm\big] \; \stackrel{(i)}{\geq} \; \left( \big(1 -
3\epsilon\big) \; (1 - \frac{\fdr}{\epsilon} \right) \; \FDPm,
\end{align*}
where step (i) uses the lower bound~\eqref{EqGoodEvents}.

Recalling that $\FDR(K) \leq \fdr$ by assumption and rearranging the
inequality, we find that
\begin{align*}
\FDPm \leq \frac{1}{ \left(1 - \frac{\fdr}{\epsilon} \right) \cdot
  \big(1 - 3\epsilon\big)} \cdot \fdr \leq \frac{2}{1 - 3\epsilon}
\cdot \fdr,
\end{align*}
where we have also used the assumed inequality $\frac{\fdr}{\epsilon}
\leq \frac{1}{2}$.  This establishes the first inequality in
line~\eqref{EqnLineOne}.

We now prove the second inequality in line~\eqref{EqnLineOne}.
Observe that the number of non-discovered signals in the top $\kstarm$
statistics can be lower bounded as
\begin{align*}
m \cdot \FNP(\kstarm) & = m - \big(\kstarm - L(\kstarm)\big) \; \geq
\; m - \big(\kstarm - \lstarm\big) \; = \; m \cdot \FNPm.
\end{align*}
Next note that conditioned on the event $\goodeventm$, we have
\begin{align*}
m \cdot \FNP\big(\kstarm\big) & \leq m \cdot \FNP\big(K\big) + \big(K
- \kstarm\big) \\
& \leq m \cdot \FNP\big(K\big) + \big(\kstarp - \kstarm\big) \\
& \leq m \cdot \bigg[\FNP\big(K\big) + 2\epsinv\big(\fdr +
  \fnr\big)\bigg] ,
\end{align*}
where we have used the fact that
\begin{align*}
\kstarp - \kstarm \leq \frac{\epsinv \big(\fdr + \fnr \big)}{1 -
  \epsinv \fdr} \leq 2\epsinv\big(\fdr + \fnr\big).
\end{align*}
Once again taking conditional expectations, dividing through by
$\P[\goodeventm]$, and using the bound on $\FNR(K)$, we find
\begin{align*}
m \cdot \FNP\big(\kstarm\big) & \leq m \cdot \bigg[\frac{\fnr}{1 - 3
    \epsilon} + 2\epsinv\big(\fdr + \fnr\big)\bigg]
\end{align*}
Putting it all together, we conclude that
\begin{align*}
  \FNPm & \leq \frac{\fnr}{1 - 3\epsilon} + 2\epsinv\big(\fdr + \fnr\big)
  \\
& \leq \bigg(\frac{1}{1 - 3\epsilon} + 2\epsinv\bigg) \cdot \big(\fdr +
  \fnr\big) \\
& \leq \bigg(\frac{2}{1 - 3\epsilon} + 4\epsinv\bigg) \cdot
  \max\big\lbrace \fdr, ~\fnr\big\rbrace .
\end{align*}


\paragraph{Argument for $+$ proxies:}  

The argument for the bounds in line~\eqref{EqnLineTwo} based on
positive-subscripted proxies is similar, so that we merely sketch
it. By reasoning similar to that used for FNR above, we can show that
\begin{align*}
\lstarp & \leq m \cdot \bigg[\frac{2}{1 - 3\epsilon} \cdot \fdr + 2
  \epsinv\big(\fdr + \fnr \big) \bigg] .
\end{align*}
Rearranging as before yields the inequality $\FDPp \leq
\bigg(\frac{4}{1 - 3\epsilon} + 4\epsinv\bigg) \cdot \max\big\lbrace
\fdr,~\fnr\big\rbrace$.  Conditioned on $\goodeventp$, we have
\begin{align*}
& m - \kstarp + \lstarp \leq m \cdot \FNP(\kstarp) \leq m \cdot
  \FNP(K) \\ &~~~~~~~~~ \implies \FNP(K) \geq \FNPp .
\end{align*}
Taking conditional expectations and dividing by the probability
$\P[\goodeventp]$, we conclude that
\begin{align*}
\FNPp & \geq \frac{1}{1 - 3\epsilon} \cdot \fnr,
\end{align*}
which completes the proof of the first inequality.  The proof of the
second inequality is analogous to the negative-subscripted case.


\subsubsection{Proof of Lemma~\ref{LemBand}}
\label{SecLemBand}
It suffices to establish the inequalities $\P [K \geq \kstarp] \leq
\epsilon$ and $\P [K \leq \kstarm ] \leq \epsilon$.  Beginning with
the first inequality, note that $\FDP\big(K\big) \geq \frac{K - m}{K}
= 1 - \frac{m}{K}$ and that the lower bound is an increasing function
of $K$. Thus, we have the inclusions
\begin{align*}
\eventwrap{K \geq \kstarp} & \subset \eventwrap{\FDP\big(K\big) \geq 1
  - \frac{m}{\kstarp}} \\
& = \eventwrap{\FDP\big(K\big) \geq \EPSINV \fdr} \\
& \subset \eventwrap{\FDP\big(K\big) \geq \EPSINV \FDR\big(K\big)}.
\end{align*}
Given this set inclusion, we have
\begin{align*}
\mprob \big[ K \geq \kstarp \big] & \leq \mprob \left[ \FDP(K) \geq
  \frac{\FDR(K)}{\epsilon} \right] \; \leq \; \epsilon,
\end{align*}
where the final line follows by Markov's inequality.

As for the second inequality, note that $\FNP\big(K\big) \geq \frac{m
  - K}{m} = 1 - \frac{K}{m}$ and that this lower bound is a decreasing
function of $K$. Thus, we have the inclusions
\begin{align*}
\eventwrap{K \leq \kstarm(\fnr)} & \subset \eventwrap{\FNP\big(K\big)
  \geq 1 - \frac{\kstarm}{m}} \\ & = \eventwrap{\FNP\big(K\big) \geq
  \EPSINV \fnr} \\ & \subset \eventwrap{\FNP\big(K\big) \geq \EPSINV
  \FNR\big(K\big)} .
\end{align*}
As before, applying Markov's inequality yields the claim.


\subsubsection{Proof of Lemma~\ref{LemOrdStatsRelation}}
\label{SecLemOrdStatsRelation}

Suppose that $W_{\nulls,(\ell)} > \transf\big(W_{\sigs,(k - \ell +
  1)}\big)$, or equivalently, $X_{\nulls,(\ell)} > X_{\sigs,(k - \ell
  + 1)}$. Define the set \mbox{$I = \big\lbrace i \colon X_{i} \geq
  X_{\nulls,(\ell)}\big\rbrace$,} and note that if $\big|I\big| \leq
k$, then necessarily $X_{\nulls,(\ell)}$ is one of the top $k$
statistics, so that $L(k) \geq \ell$. But, by the hypothesis and the
definition of order statistics,
\begin{align*}
\big|I \cap S\big| \leq k - \ell ~~ \text{and} ~~ \big|I \cap
\nulls\big| = \ell.
\end{align*}
Thus $\big|I\big| = \big|I \cap S\big| + \big|I \cap \nulls\big| \leq
k$, as required.

We now turn to the converse implication.  Concretely, fixing some $k
\in [n] \defn \{1, 2, \ldots, n \}$ for which $L(k) \geq \ell$, we
prove that $W_{\nulls,(\ell)} > \transf\big(W_{\sigs,(k - \ell +
  1)}\big)$.  Let $i_{k} \in [n]$ be the $k^{th}$-largest rank
statistic---that is, the index corresponding to the order statistic
$X_{(k)}$.  and we break our analysis into two cases, depending on
whether $i_k \in \sigs$ or $i \in \nulls$.

\paragraph{Case 1, $i_k \in \sigs$:}
In this case, since there are at most $k - \ell$ signals in the top
$k$ statistics, we must have $X_{i_k} \geq X_{\sigs,(k - \ell)}$. On
the other hand, for any $j \in \nulls$ such that $X_{j}$ falls in the
top $k$, we must have $X_{j} > X_{i_k}$. Since there are at least
$\ell$ such indices, we conclude
\begin{align*}
W_{\nulls,(\ell)} = X_{\nulls,(\ell)} > X_{i_k} \geq X_{\sigs,(k -
  \ell)} = \transf\big(W_{\sigs,(k - \ell)}\big) >
\transf\big(W_{\sigs,(k - \ell + 1)}\big) .
\end{align*}
Rearranging yields the claim.


\paragraph{Case 2, $i_k \in \nulls$:}
Since the number of signals in the top $k$ is $< k - \ell + 1$, it
must be that $X_{\sigs,(k - \ell + 1)} < X_{i_{k}} \leq
X_{\nulls,(\ell)}$. Rearranging again gives the claim.


\subsection{Proof of Corollary~\ref{ThmIIDNormal}}
\label{SecIIDNormalProof}

Let $\gaussord_{k,n}$ denote the expected value of the
$k^{\text{th}}$-largest value in a sample of $n$ independent standard
Gaussians.  Recalling the definition~\eqref{EqnConcentration} of the
concentration function, classical results on Gaussian order statistics
ensure that we can apply the concentration assumption with
\mbox{$\ordD_{k,n}(\epsilon) \defn \sqrt{2 \log \frac{2}{\epsilon}}$.}
Although this specification is not the sharpest possible, it suffices
for our purposes.

Our proof of Corollary~\ref{ThmIIDNormal} is based primarily on
comparing Gaussian order statistics to $\mu$. In particular, we wish
to establish that the inequality
\begin{align}
\label{EqIIDNormalMain}
\gaussord_{\lstar + 1,n - m} + \gaussord_{m - \kstar + \lstar, m} >
\mu - 2\sqrt{2 \log \frac{6}{\epsilon}}
\end{align}
holds with the choices $(\lstar,~\kstar) =
(\lstarm(\fnr),~\kstarm(\fnr))$, or $(\lstar,~\kstar) =
(\lstarp(\fdr),~\kstarp(\fdr))$.  The proof is identical for these two
cases, so we simply use the shorter $(\lstar, \kstar)$ notation
throughout.

Taking inequality~\eqref{EqIIDNormalMain} as given for the moment, we
first use it to prove Corollary~\ref{ThmIIDNormal}.  In order to do
so, we require the following:
\begin{lem}
\label{LemNormalOrdStats}
We have
\begin{align}
\frac{\sqrt{2 \log \frac{n}{k}}}{\xi_{k,n}} & = 1 \pm o(1) \quad
\mbox{for $k = 1, \ldots, 2, \ldots, n$,}
\end{align}
where the $o(1)$ decay holds as $n \rightarrow \infty$ and/or $k
\rightarrow \infty$.
\end{lem}

We now proceed with the proof, suppressing $n$ subscripts throughout
so as to avoid clutter.  Combining inequality~\eqref{EqIIDNormalMain},
Lemma~\ref{LemNormalOrdStats}, and the fact that $m \ll n$, we find
that
\begin{align*}
\mu > \gaussord_{\lstar + 1,n - m} + \gaussord_{\kstar - \lstar, m} -
c_{1} \geq \big(1 - o(1) \big)\bigg[\sqrt{\log{\frac{n}{\lstar + 1}}}
  + \sqrt{\log{\frac{m}{m - \kstar + \lstar}}}\bigg],
\end{align*}
where $c_1 > 0$ is a constant that may depend on $\epsilon$.
\footnote{To clarify a subtle point that we have elided:
  Lemma~\ref{LemNormalOrdStats} requires that $k < \frac{n}{\log{n}}$,
  so we need to check $\lstar + 1 \leq \frac{n}{\log{n}}$ and $m -
  \kstar + \lstar \leq \frac{n}{\log{n}}$. Since $\max \big\lbrace
  \lstar,~m - \kstar + \lstar\big\rbrace \leq \max \big\lbrace
  q,~\fnr\big\rbrace \cdot m$ by Theorem~\ref{ThmProxies}, this
  condition is in fact easily verified under the given scalings for
  $\fdr$ and $\fnr$.}

We now invoke Theorem~\ref{ThmProxies} with $(\lstar, \kstar) =
(\lstarm, \kstarm)$ or $(\lstar, \kstar) = (\lstarp, \kstarp)$,
according to whether we are in the regime $\fnr \geq \fdr$ or vice
versa.  Applying the theorem, rearranging, and substituting the value
of $\mu$ yields
\begin{align}
  \label{EqIIDNormalSecondary}
\sqrt{\log{\frac{n}{\cproxy \fdr m + 1}}} +
\sqrt{\log{\frac{m}{\cproxy \fnr m}}} \leq \big(1 - o(1) \big)^{-1}
\sqrt{2r \log{n}}  \; = \big(1 + o(1) \big) \sqrt{2 r \log n}.
\end{align}

We claim that it suffices to prove that $\min\big\lbrace
\fdr,~\fnr\big\rbrace \cdot m \rightarrow \infty$.  Indeed, under this
scaling, for large enough $(n,m)$, we would ahve
\begin{align*}
\sqrt{2 \log \frac{n}{2\cproxy \fdr m}} + \sqrt{2
  \log\frac{m}{2\cproxy \fnr m}} \leq (1 + o(1)) \cdot \sqrt{2r
  \log{n}} .
\end{align*}
Substituting the assumed scalings $\fdr = c n^{-\kappa_\fdr}$ and
$\fnr = c n^{-\kappa_\fnr}$ then yields
\begin{align*}
\sqrt{2\big(\betapar + \kappa_\fdr\big) \log{n} +
  \log\frac{1}{2\cproxy}} + \sqrt{2\kappa_\fnr \log{n} +
  \log\frac{1}{2\cproxy}} \leq \big(1 + o(1) \big) \cdot \sqrt{2r
  \log{n}},
\end{align*}
and letting $n \rightarrow \infty$ yields the claimed inequality
$\sqrt{\betapar + \kappa_\fdr} + \sqrt{\kappa_\fnr} \leq \sqrt{r}$.

It remains to prove that $\min \{ \fdr, ~\fnr \} \cdot m \rightarrow
\infty$, and we split our analysis into two cases.

\paragraph{Case 1:}  Suppose first that $\fdr \leq \fnr$ and assume by way of
contradiction that there exists a constant $c_{2}$ such that $\fdr m
\leq c_{2}$ for all $n$.  Combined with the
inequality~\eqref{EqIIDNormalSecondary}, we find that
\begin{align*}
\sqrt{2 \log \frac{n}{\cproxy c_{2} + 1}} + \sqrt{2
  \log\frac{m}{\cproxy \fnr m}} \leq \big(1 + o(1) \big) \cdot
\sqrt{2r \log{n}} .
\end{align*}
Since $r < 1$, this inequality cannot hold once $n$ is large enough,
which establishes the desired contradiction.


\paragraph{Case 2:}  Turning to the other case, suppose that
$\fnr \leq \fdr$, and assume by way of contradiction that there exists
a constant $c_{2}$ such that $\fnr m \leq c_{2}$ for all $n$. In this
case, again by inequality~\eqref{EqIIDNormalSecondary}, we have
\begin{align*}
\sqrt{2 \log \frac{n}{\cproxy \fdr m + 1}} + \sqrt{2
  \log\frac{m}{\cproxy c_{2}}} \leq \big(1 + o(1) \big) \cdot \sqrt{2r
  \log{n}} .
\end{align*}
On the other hand, for a suitable choice of $c_{3}$, we have
\begin{align*}
\sqrt{2 \log \frac{n}{\cproxy \fdr m + 1}} + \sqrt{2
  \log\frac{m}{\cproxy c_{2}}} & \geq \sqrt{2 \log \frac{n}{4\cproxy
    m}} + \sqrt{2 \log\frac{m}{\cproxy c_{2}}} \\
& \geq \sqrt{2 \betapar \log{n} +
  \log\frac{1}{c_{3}}} + \sqrt{2\big(1 - \betapar \big) \log{n} +
  \log\frac{1}{c_{3}}} .
\end{align*}
Putting together the pieces, we have shown that
\begin{align*}
\sqrt{2\betapar \log{n} + \log\frac{1}{c_{3}}} + \sqrt{2\big(1 -
  \betapar \big) \log{n} + \log\frac{1}{c_{3}}} & \leq \big(1 + o(1)
\big) \cdot \sqrt{2r \log{n}} .
\end{align*}
Since $\sqrt{\betapar} + \sqrt{1 - \betapar} > 1 > \sqrt{r}$, this
inequality cannot hold once $n$ is sufficiently large, which establishes
the desired contradiction in this case.


\subsubsection{Proof of inequality~\eqref{EqIIDNormalMain}}

Applying Lemma~\ref{LemDerandom} with $\ordD_{\nulls,k} =
\ordD_{\sigs,k} = \sqrt{2\log\frac{2}{\epsilon}}$, we find that
\begin{align*}
\mu + \gaussord_{\kstar - \lstar,m} + \sqrt{2 \log\frac{2}{\epsilon}}
> \gaussord_{\lstar + 1, n -m} - \sqrt{2 \log\frac{2}{\epsilon}} .
\end{align*}
Rearranging yields
\begin{align*}
\mu > \gaussord_{\lstar + 1, n -m} - \gaussord_{\kstar - \lstar,m} - 2
\sqrt{2 \log\frac{2}{\epsilon}} .
\end{align*}
Since the Gaussian distribution is symmetric around zero, we can
replace $-\gaussord_{\kstar - \lstar,m}$ by $\gaussord_{m - \kstar +
  \lstar, m}$, which yields the desired inequality.


\subsection{Proof of Corollary~\ref{ThmIIDNormalScale}}

By analogy to the notation in Section~\ref{SecIIDNormalProof}, let
$\chi_{k,n}$ denote the expected value of the
\mbox{$k^{\text{th}}$-largest} value in a sample of $n$ variables,
each of which is the absolute value of a standard Gaussian.  Other
notational conventions are also preserved. In particular, we let
$(\kstar,\lstar)$ stand in for either $(\kstarm,~\lstarm)$ or
$(\kstarp,~\lstarp)$, depending on whether $\fnr \geq \fdr$ or vice
versa. We also suppress $n$ subscripts.

By applying Lemma~\ref{LemDerandom} in this case, we find that
\begin{align*}
\sigma \cdot \gaussord_{\kstar - \lstar,m} + \sqrt{2 \log
  \frac{6}{\epsilon}} \geq \gaussord_{\lstar + 1,n - m} - \sqrt{2 \log
  \frac{6}{\epsilon}},
\end{align*}
and rearranging yields
\begin{align}
\label{EqNormalScaleMain-2}
\sigma \geq \gaussord_{\kstar - \lstar,m}^{-1}\bigg[\gaussord_{\lstar
    + 1,n-m} + 2\sqrt{2 \log\frac{6}{\epsilon}}\bigg] .
\end{align}

We now need to evaluate the $\gaussord$ values, and we make use of the
following result due to Gordon et al.~\cite{Gor07Interpolated}:
\begin{lem}
\label{LemNormalScaleOrdStats}
For all $k \geq n/2$, we have $\gaussord_{k,n} \leq \sqrt{2\pi} \cdot
\frac{n - k + 1}{n + 1}$.  Moreover, we have
\begin{align}
 \frac{\sqrt{2 \log \frac{n}{k}}}{\gaussord_{k,n}} & = 1 \pm o(1)
 \qquad \mbox{for all $k = 1, 2, \ldots, \lfloor \frac{n}{\log n}
   \rfloor$,}
\end{align}
where the $o(1)$ scaling holds as $n$ and possibly $k$ go to infinity.
\end{lem}

Suppose for now that $\lstar + 1 \leq \frac{n - m}{\log{n - m}}$ and
that $\kstar - \lstar \geq \frac{m}{2}$. Recall from
Theorem~\ref{ThmProxies}, applied with the appropriate choice of $+$
or $-$ proxies, that
\begin{align}
\label{EqNormalScaleReminder}
\lstar \leq \cproxy \fdr m ~~ \text{and} ~~ m - \kstar + \lstar \leq
\cproxy \fnr m .
\end{align}
Consequently, we have the lower bound
\begin{subequations}
\begin{align}
\gaussord_{\lstar + 1,n-m} & \geq \big(1 - o(1) \big) \cdot
\sqrt{2\log \frac{n - m}{\lstar + 1}} \nonumber \\
& \geq \big(1 - o(1) \big) \cdot \sqrt{2\log \frac{\frac{n}{m} -
    1}{\cproxy \fdr + 1/m}} \nonumber \\
\label{EqnRameshA}
& \geq (1 - o(1) \big) \cdot \sqrt{2 \betapar \log{n} +
  2\log\bigg(\fdr + \frac{1}{m}\bigg)^{-1}}.
\end{align}
On the other hand, Lemma~\ref{LemNormalScaleOrdStats} also implies
that
\begin{align}
\gaussord_{\kstar - \lstar,m} & \leq \sqrt{2\pi} \cdot \frac{m -
  \kstar + \lstar + 1}{m + 1} \nonumber \\
& \leq \sqrt{2\pi} \cdot\bigg( \frac{m - \kstar - \lstar}{m} +
\frac{1}{m}\bigg) \nonumber \\
\label{EqnRameshB}
& \leq \sqrt{2\pi}\cproxy \cdot \big(\fnr +\frac{1}{m}\big) .
\end{align}
\end{subequations}
Combining the bounds~\eqref{EqnRameshA} and~\eqref{EqnRameshB}
with inequality~\eqref{EqNormalScaleMain-2}, we find that
\begin{align*}
\sigma \geq \frac{1}{\sqrt{2\pi}\cproxy} \cdot \big(1 - o(1)
\big)\bigg(\fnr + \frac{1}{m}\bigg)^{-1}\bigg[\sqrt{2\betapar \log{n}
    + 2\log\bigg(\fdr + \frac{1}{m}\bigg)^{-1}} + 2\sqrt{2
    \log\frac{6}{\epsilon}}\bigg] .
\end{align*}
Since $\betapar \geq \rho > 0$, we have that $\sqrt{2\beta \log{n} +
  2\log\bigg(\fdr+ \frac{1}{m}\bigg)^{-1}} \gg 2\sqrt{2
  \log\frac{6}{\epsilon}}$, so that
\begin{align*}
\sigma \geq \frac{1}{\sqrt{2\pi}\cproxy} \cdot \big(1 - o(1) \big)
\bigg(\fnr + \frac{1}{m}\bigg)^{-1}\sqrt{2\betapar \log{n} + 2 \log
  \bigg(\fdr + \frac{1}{m}\bigg)^{-1}},
\end{align*}
as claimed.  

Under the specified scalings, we have $\fdr + \frac{1}{m} \approx
n^{\kappa_\fdr}$ and $\fnr + \frac{1}{m} \approx n^{\kappa_\fnr}$,
which directly implies the comparison
\begin{align*}
\sigma \gtrsim n^{\kappa_\fnr}\sqrt{2\big(\betapar + \kappa_\fdr \big)
  \log{n}} .
\end{align*}

We now need to verify that $\lstar + 1 \leq \frac{n - m}{\log\big(n -
  m\big)}$ and $\kstar - \lstar \geq \frac{m}{2}$. From the
inequalities~\eqref{EqNormalScaleReminder}, we deduce that
\begin{align*}
\lstar + 1 \leq \cproxy \fdr m + 1 ~~ \text{and} ~~ \kstar - \lstar
\geq \big(1 - \cproxy \fnr \big) m .
\end{align*}
Note that, by the assumption that $\betapar \geq \rho$, we have $m
\leq n^{1 - \rho}$, we have $\frac{n - m}{\log\big(n - m\big)} \geq
\frac{n}{2 \log{n}}$ for large enough $n$.  On the other hand, we also
have $\fdr m + 1 \leq n^{1 - \rho} + 1 \leq \frac{n}{2\log{n}}$ (once
$n$ large enough---say for all $n$ such that $n^{\rho} \geq 4
\log{n}$, for instance). For the second case, recall the assumption
$\cproxy \fnr \leq \frac{1}{3} < \frac{1}{2}$, from which the claim
follows.


\subsection{Proof of Corollary~\ref{ThmSpikedNormal}}

At a high level, this proof involves reducing to an independent model
with altered variances and using Lipschitz concentration to verify the
closeness condition of Lemma~\ref{LemTransfer}.

We carry out the reduction in two steps: in Step 1, we reduce to a
model with dependence \emph{only} between nulls and signals, and then
in Step 2, we reduce to an independent model.  The models in Steps 1
and 2 are Gaussian models with covariance matrices $\Sigma'$ and
$\Sigma''$, respectively, of the form
\begin{align*}
\Sigma'_{ij} = \begin{cases} 1 - \rhonull + \rhocross &~\text{if}~i
  = j,~i \in \nulls, \\ 1 - \rhosig + \rhocross &~\text{if}~i = j,~i
  \notin \nulls, \\ \rhocross &~\text{if}~i \in \nulls,~j \notin
  \nulls, \\ 0 &~\text{o.w.}
\end{cases}
\quad \mbox{and} \quad \Sigma''_{ij} & = \begin{cases} 1 - \rhonull
  &~\text{if}~i = j,~i \in \nulls, \\ 1 - \rhosig &~\text{if}~i = j,~i
  \notin \nulls, \\ 0 &~\text{o.w.}
\end{cases}
\end{align*}
We let $W'$ and $W''$ corresponding the corresponding Gaussian random
vectors in $\R^{n}$. The shifts associated with these models are set
to be constant-scale perturbations of $\mu$, so that overall, we have
the two models
\begin{align*}
\model' = \left(\PW',~\mu + c_{1}\sqrt{2 \log\frac{c_{2}}{\epsilon}}
\right) ~~ \text{and} ~~ \model'' = \left (\PW'',~\mu +
2c_{1}\sqrt{2\log\frac{c_{2}}{\epsilon}} \right),
\end{align*}
where $\PW'$ and $\PW''$ are the Gaussian distributions associated
with $\Sigma'$ and $\Sigma''$.  We also introduce the convenient
shorthand notation
\begin{align*}
{\lstar}' = \lstar \left( \model',~\frac{\epsilon}{3} \right) ~~
\text{and} ~~ {\lstar}'' = \lstar \left(\model'',~\frac{\epsilon}{9}
\right) .
\end{align*}

The main idea of the proof is to represent the $W$ variables as
functions of higher-dimensional Gaussians. This representation is
helpful in decoupling the test statistics from each other.  Basically,
the constant covariance within the nulls and signals, and across the
two, allows us to represent each test statistic as independent of all
the others after conditioning on three standard Gaussians that contain
all the shared randomness: one each for the within-nulls,
within-signals, and between-nulls-and-signals randomness.  More
precisely, we can write
\begin{align*}
W_{i} & = \begin{cases} \sqrt{1 - \rhonull} \cdot U_{i} +
  \sqrt{\rhonull - \rhocross} \cdot V_{0} + \sqrt{\rhocross} \cdot
  \Vcross &~\text{if}~ i \in \nulls, \\ \sqrt{1 - \rhosig} \cdot U_{i}
  + \sqrt{\rhosig - \rhocross} \cdot V_{1} + \sqrt{\rhocross} \cdot
  \Vcross &~\text{if}~i \notin \nulls
\end{cases}, \\
W_{i}' & = \begin{cases} \sqrt{1 - \rhonull} \cdot U_{i}' +
  \sqrt{\rhocross} \cdot \Vcross' &~\text{if}~ i \in \nulls,
  \\ \sqrt{1 - \rhosig} \cdot U_{i}' + \sqrt{\rhocross} \cdot \Vcross'
  &~\text{if}~i \notin \nulls
\end{cases}, \\
W_{i}'' & = \begin{cases} \sqrt{1 - \rhonull} \cdot U_{i}''
  &~\text{if}~ i \in \nulls, \\ \sqrt{1 - \rhosig} \cdot U_{i}''
  &~\text{if}~i \notin \nulls
\end{cases},
\end{align*}
The link functions that connect the $U$ and $V$ variables to the order
statistics of the $W$ variables are given, in the three cases, by
\begin{align*}
a_{0,\ell} \big(u,~v_{0:1},~v_{c}\big) & = \big(\sqrt{1 - \rhonull}
\cdot u + \sqrt{\rhonull - \rhocross} \cdot v_{0} + \sqrt{\rhocross}
\cdot v_{c}\big)_{(\ell)}, \\
b_{0,\ell} \big(u,~v_{c}\big) & = \big(\sqrt{1 - \rhonull} \cdot u +
\sqrt{\rhocross} \cdot v_{c} \big)_{(\ell)}, \\
c_{0,\ell}\big(u\big) & = \sqrt{1 - \rhonull} \cdot u_{(\ell)} ,
\end{align*}
and similarly for the signals, for which we denote the link functions
by $a_{1}$, $b_{1}$, and $c_{1}$.

Our first aim is to prove $\lstar \geq {\lstar}'$ using
Lemma~\ref{LemTransfer}.   We begin by observing that
\begin{align*}
\bigg|a_{0,\ell}\big(u,~v_{0:1},~\vcross\big) -
b_{0,\ell}\big(u,~\vcross\big)\bigg| \leq \big|v_{0}\big| ,
\end{align*}
from which it follows that
\begin{align*}
\eventwrap{\bigg|W_{\nulls,({\lstar}')} -
  W_{\nulls,({\lstar}')}'\bigg| \leq \Delta} & =
\eventwrap{\bigg|a_{0,{\lstar}'}\big(U,~V_{0:1},~V_{c}\big) -
  b_{0,{\lstar}'}\big(U',~V_{c}'\big)\bigg| \leq \Delta} \\
& \supset
\eventwrap{|V_{0}| \leq \frac{\Delta}{2}} \cap
\eventwrap{\bigg|b_{0,{\lstar}'}\big(U,~V_{c}\big) -
  b_{0,{\lstar}'}\big(U',~V_{c}'\big)\bigg| \leq \frac{\Delta}{2}} \\
& \supset \eventwrap{|V_{0}| \leq \frac{\Delta}{2}} \cap
\eventwrap{\bigg|b_{0,{\lstar}'}\big(U,~V_{c}\big) -
  \E\big[b_{0,{\lstar}'}\big(U,~V_{c}\big)\big]\bigg| \leq
  \frac{\Delta}{4}} \\
& \hspace{1.05in} \cap
\eventwrap{\bigg|b_{0,{\lstar}'}\big(U',~V_{c}'\big) -
  \E\big[b_{0,{\lstar}'}\big(U',~V_{c}'\big)\big]\bigg| \leq
  \frac{\Delta}{4}}.
\end{align*}
A similar analysis yields that the event
$\eventwrap{\bigg|W_{\sigs,({\kstar}' - {\lstar}')} -
  W_{\sigs,({\kstar}' - {\lstar}')}'\bigg| \leq \Delta}$ contains the
event
\begin{align*}
\eventwrap{|V_{1}| \leq \frac{\Delta}{2}} \cap
\eventwrap{\bigg|b_{1,{\kstar}' - {\lstar}'}\big(U,~V_{c}\big) -
  \E\big[b_{1,{\kstar}' - {\lstar}'}\big(U,~V_{c}\big)\big]\bigg| \leq
  \frac{\Delta}{4}} \\
\qquad \qquad \cap \eventwrap{\bigg|b_{1,{\kstar}' -
    {\lstar}'}\big(U',~V_{c}'\big) - \E\big[b_{1,{\kstar}' -
      {\lstar}'}\big(U',~V_{c}'\big)\big]\bigg| \leq \frac{\Delta}{4}}.
\end{align*}
By Lipschitz concentration, we may choose $\Delta = c_{1}'\sqrt{2 \log
  \frac{c_{2}'}{\epsilon}}$ such that
\begin{align*}
\max\Bigg\lbrace \P\bigg(\big|W_{\nulls,({\lstar}')} -
W_{\nulls,({\lstar}')}'\big| > \Delta\bigg),~
\P\bigg(\big|W_{\sigs,({\kstar}' - {\lstar}')} - W_{\sigs,({\kstar}' -
  {\lstar}')}'\big| > \Delta\bigg)\Bigg\rbrace \leq \frac{\epsilon}{3}.
\end{align*}
In other words, $\model$ and $\model'$ are both $\big(\Delta,~\Delta,~
\lstarm,~\frac{\epsilon}{3}\big)$ and $\big(\Delta,~\Delta,~
\lstarp,~\frac{\epsilon}{3}\big)$ close.  Consequently, if we choose
$c_{1} \geq 2c_{1}'$ and $c_{2} \geq c_{2}'$, then applying
Lemma~\ref{LemTransfer} guarantees that ${\lstar}' \geq \lstar$.

Next observe that $\big | b_{0,\ell}(u,~v_{c}) - c_{0,\ell}(u) \big |
\leq |v_{c}|$.  We may therefore apply a variant of the previous
argument with $\Delta' = c_{1}''\sqrt{2 \log
  \frac{c_{2}''}{\epsilon}}$ to show
\begin{align*}
\max \Bigg \lbrace \P\bigg(\big|W_{\nulls,({\lstar}'')}' -
W_{\nulls,({\lstar}'')}''\big| > \Delta'\bigg),~
\P\bigg(\big|W_{\sigs,(({\kstar}'' - {\lstar}'')}' -
W_{\sigs,({\kstar}'' - {\lstar}'')}''\big| > \Delta'\bigg)\Bigg\rbrace
\leq \frac{\epsilon}{9} .
\end{align*}
We then find by Lemma~\ref{LemTransfer} that ${\lstar}'' \geq
{\lstar}'$ provided that $c_{1} \geq 2c_{1}''$ and $c_{2} \geq
c_{2}''$.

Combining the two pieces of our argument, we are guaranteed to have
$\lstar \geq {\lstar}''$ as long $c_{1} \geq 2 \cdot \max\big\lbrace
c_{1}',~c_{1}''\big\rbrace$ and $c_{2} \geq \max\big\lbrace
c_{2}',~c_{2}''\big\rbrace$.  Since we have now reduced to the
independent case with $\mu$ changed by a constant, applying suitably
rescaled version of Theorem~\ref{ThmIIDNormal} yields the conclusion
of the corollary.


\subsection{Proof of Corollary~\ref{ThmGroupedNormal}}

By our previous arguments for order statistics of Gaussians, we know
that for the grouped Gaussian model, Case I of the concentration
assumption~\eqref{EqnConcentration} holds with
\begin{align*}
  \ordD_{\nulls,\lstar} = \ordD_{\sigs,\kstar - \lstar + 1} =
  c_{1}\sqrt{2\log\frac{c_{2}}{\epsilon}}
\end{align*}
If we choose the constants $c_1, c_2$ sufficiently large, we can
ensure that Case II of the concentration
condition~\eqref{EqnConcentration} holds at the same scale in a
modified form $\model'$ of the model in which the dependence between
nulls and signals is broken and the shift is altered to $\mu -
c_{1}'\sqrt{2 \log \frac{c_{2}'}{\epsilon}}$.  We may therefore apply
Lemma~\ref{LemDecouple} to obtain that $\lstar \geq {\lstar}'$, where
${\lstar}'$ is computed in $\model'$.

Since Case II of the concentration condition~\eqref{EqnConcentration}
holds in the new model, we may apply
Lemma~\ref{LemDerandom}. Specifically, if we set $T_{0} =
\nulls^{(0)}$ and $T_{1} = \nulls\setminus T_{0}$, then
\begin{align*}
\mu - c_{1}'\sqrt{2 \log \frac{c_{2}'}{\epsilon}} & \geq
\E\big[W_{\nulls,({\lstar}' + 1)}'\big] + \E\big[W_{\sigs,(m - \kstar
    + {\lstar}' + 1)}'\big] - 2c_{1}'\sqrt{2 \log
  \frac{c_{2}'}{\epsilon}} \\
& \geq \E\big[W_{T_{0},({\lstar}' + 1)}'\big] + \E\big[W_{\sigs,(m -
    \kstar + {\lstar}' + 1)}'\big] - 2c_{1}'\sqrt{2 \log
  \frac{c_{2}'}{\epsilon}} \\
& = \gaussord_{{\lstar}' + 1,|T_{0}|} + \gaussord_{m - \kstar +
  {\lstar}' + 1, m} - 2c_{1}'\sqrt{2 \log \frac{c_{2}'}{\epsilon}} .
\end{align*}
Since $\etapar < \betapar$, we have $|T_{0}| \geq n -
(\frac{1}{\lehmann} + 1) m \geq n - (n^{\etapar} + 1)n^{1 - \betapar}
\geq n - 2n^{1 - \etapar } = \big(1 - o(1)\big)n$, and an application
of Lemma~\ref{LemNormalOrdStats} yields the claim.


\subsection{Proof of Corollary~\ref{ThmLehmannAlt}}

In order to simplify the proof, it is convenient to pass to an
equivalent model.  Consider the new random vector $V = (V_1, \ldots,
V_n)$ with components $V_{i} \defn \log \frac{1}{1 - W_{i}}$.  Note
that $V_{i}$ is distributed as a standard exponential and that if we
define the transformation function $g(v) = v/\lehmann= Av$ with $A
\colon = 1/\lehmann$, then
\begin{align*}
\log \frac{1}{1 - f(W_{i})} = g(V_{i}) .
\end{align*}
With this set-up, the test statistics in the new model are related to
the test statistics in the original model by the transformation $x
\mapsto \log \frac{1}{1 - x}$. Since this transformation is monotonic,
any \topk procedure for one can be translated into a \topk procedure
for the other, with no change in performance.  Likewise, the proxy
values are the same for all $\fdr$ and $\fnr$.  In summary, the two
models are equivalent for our purposes.

As in previous proofs, we use $(\kstar,~\lstar)$ as a stand-in for
$(\kstarm,~\lstarm)$ or $(\kstarp,~\lstarp)$, and we suppress $n$
subscripts.  We claim that it is sufficient to show that
\begin{subequations}
\begin{align}
\label{EqVMinus}
\P\bigg[ V_{\nulls,(\lstar + 1)} \leq v_{-}\bigg] \leq
\frac{\epsilon}{3}, \\
\label{EqVPlus} 
\P\bigg[ g(V_{\sigs,(\kstar - \lstar)}) \geq v_{+}\bigg] \leq
\frac{\epsilon}{3},
\end{align}
\end{subequations}
where
\begin{align*}
v_{+} = \frac{1}{\lehmann} \frac{\etapar}{1 - \etapar} ~~ \text{and} ~~
v_{-} = \log\Bigg(\frac{1}{\cproxy\pi_{1} \fdr}\bigg(1 + 4
\log{\frac{3}{\epsilon}}\bigg)^{-1}\Bigg).
\end{align*}

Taking inequalities~\eqref{EqVMinus} and~\eqref{EqVPlus} as given for
the moment, by the definition of $\lstar$, we have
\begin{align*}
\P\bigg [ g(V_{\sigs,(\kstar - \lstar)}) > V_{\nulls,(\lstar + 1)}
  \bigg] \geq \epsilon .
\end{align*}
On the other hand, combining the two bounds above, we see that
\begin{align*}
\P\bigg(V_{\nulls,(\lstar + 1)} \leq v_{-}\bigg) +
\P\bigg(g(V_{\sigs,(\kstar - \lstar)}) \geq v_{+}\bigg) \leq
\frac{2\epsilon}{3} < \epsilon .
\end{align*}
It follows that
\begin{align*}
\P\bigg [ v_{+} > g(V_{\sigs,(\kstar - \lstar)}) > V_{\nulls,(\lstar +
    1)} > v_{-}\bigg] > 0,
\end{align*}
so $v_{-} \leq v_{+}$. Rearranging yields
\begin{align*}
\frac{1}{\lehmann} \geq \frac{1 - \etapar}{\etapar} \cdot
\log\Bigg(\frac{1}{\cproxy\pi_{1} \fdr}\bigg(1 + 4
\log{\frac{3}{\epsilon}}\bigg)^{-1}\Bigg),
\end{align*}
as claimed. \\

\noindent The only remaining detail is to prove
inequalities~\eqref{EqVMinus} and~\eqref{EqVPlus}.


\subsubsection{Proof of inequality~\eqref{EqVMinus}}

Applying Lemma 4.3 from the paper~\cite{Bou12Order} yields
\begin{align*}
\P \bigg[ V_{\nulls,(\lstar + 1)} \leq \log \frac{n}{\lstar} - z
  \bigg] & \leq \exp\bigg(-\frac{\lstar\big(e^{z} - 1\big)}{4}\bigg)
\qquad \mbox{for each $z > 0$.}
\end{align*}
In particular, choosing $z = \log\bigg(1 + \frac{4 \log
  \frac{3}{\epsilon}}{\lstar}\bigg) \leq \log\bigg(1 + 4
\log{\frac{3}{\epsilon}}\bigg)$, we deduce that
\begin{align*}
\P \bigg[ V_{\nulls,(\lstar + 1)} \leq
  \log\Bigg(\frac{n}{\lstar}\bigg(1 + 4 \log
  \frac{3}{\epsilon}\bigg)^{-1}\Bigg)\bigg] & \leq \epsilon/3 .
\end{align*}
We complete the proof by noting that $\frac{n}{\lstar} \geq
\frac{n}{\cproxy \fdr m} = \frac{1}{\cproxy \pi_{1} \fdr}$.


\subsubsection{Proof of inequality~\eqref{EqVPlus}}

The proof is based on the fact that $V_{i} = \log \frac{1}{U_{i}}$
where $U_{i}$ is a uniform random variable. Let $U_{\sigs,(j)}$ denote
the $j^{\text{th}}$-\emph{smallest} value in the sample, which follows
a beta distribution with parameters $j$ and $m - j + 1$. We thus have
\begin{align*}
\E \big[U_{\sigs,(\kstar - \lstar)}\big] &= \frac{\kstar - \lstar}{m +
  1}, \qquad \mbox{and} \\
\Var\big[U_{\sigs,(\kstar - \lstar)}\big] & \leq
\E\big[U_{\sigs,(\kstar - \lstar)}\big] \cdot \bigg(1 -
\E\big[U_{\sigs,(\kstar - \lstar)}\big]\bigg) \cdot \frac{1}{m}.
\end{align*}
Applying Chebyshev's inequality yields
\begin{align*}
\P\bigg(U_{\sigs,(\kstar - \lstar)} \leq \frac{\kstar - \lstar}{m + 1}
- \sqrt{\frac{3}{\epsilon} \cdot \frac{\kstar - \lstar}{m + 1} \cdot
  \frac{m - \kstar + \lstar}{m + 1} \cdot \frac{1}{m}}\bigg) \leq
\epsilon / 3.
\end{align*}
Using the fact that $\frac{m - k + \ell}{m + 1} \leq c_{0}\fnr$, we thus
have
\begin{align}
\label{eq:lehmann-unif-order-upper-bd}
\P \bigg[ U_{\sigs,(\kstar - \lstar)} \leq 1 - c_{0}\fnr - 1/m -
  \sqrt{c_{\epsilon/3}c_{0}\fnr/m}\bigg] = \P\bigg[ U_{\sigs,(\kstar -
    \lstar)} \leq 1 - \etapar \bigg] \leq \epsilon/3 .
\end{align}
We now note that
\begin{align*}
V_{\sigs,(\kstar - \lstar)} = \log \frac{1}{U_{\sigs,(\kstar -
    \lstar)}} & = \log \frac{1}{1 - (1 - U_{\sigs,(\kstar - \lstar)})}
\\
& = \log\bigg(1 + \frac{1 - U_{\sigs,(\kstar -\lstar)}}{1 - (1 -
  U_{\sigs,(\kstar - \lstar)})}\bigg) \\ & \leq \frac{1 -
  U_{\sigs,(\kstar -\lstar)}}{1 - (1 - U_{\sigs,(\kstar - \lstar)})} .
\end{align*}
Applying the bound~\eqref{eq:lehmann-unif-order-upper-bd} yields
$\P\bigg[V_{\sigs,(\kstar - \lstar)} \geq \frac{\etapar}{1 -
    \etapar}\bigg] \leq \epsilon/3$, as claimed. This inequality
completes the proof since $g$ is monotonically increasing and
$g\big(\frac{\etapar}{1 - \etapar}\big) = \frac{1}{\lehmann} \cdot
\frac{\etapar}{1 - \etapar}$.


\section{Discussion}

In this paper, we introduced a framework for establishing the
tradeoffs between the false discovery rate (FDR) and the false
non-discovery rate (FNR) in multiple testing problem. While a
substantial literature on the multiple testing problem has developed,
comparatively little has been established about the fundamental
tradeoffs between these two types of errors, or about the fundamental
limits on the combined risk (a weighted combination of FDR and
FNR). Moreover, this problem does not appear to be amenable to the
standard techniques for establishing lower bounds used in estimation
theory, for instance.

The framework we have put forward is fairly general, not being
sensitive to the analytic form of the test statistic distributions or
on the dependence structure between the test statistics. Instantiated
for models previously studied in the literature, our general results
recover and extend lower bounds previously proven using methods more
tailored to Gaussian-like models. Furthermore, the lower bounds
predicted by our theory can be numerically simulated for any given
model, an unusual feature useful for both further theoretical work and
potential applications.


\section*{Acknowledgements}

This material is based upon work supported in part by the Army
Research Office under contract/grant number W911NF-16-1-0368 to MJ.
This is part of a collaboration between US DOD, UK MOD and UK
Engineering and Physical Research Council (EPSRC) under the
Multidisciplinary University Research Initiative to MJ.  Portions of
this work were also supported by Office of Naval Research grant DOD
ONR-N00014-18-1-2640 and National Science Foundation grant
NSF-DMS-1612948 to MJW.


\bibliographystyle{plainnat} \bibliography{fdr_general}

\appendix

\section{Proofs of technical tools}
\label{AppTechnical}

In this appendix, we collect the proofs of various technical lemmas
used in the paper.


\subsection{Proof of Lemma~\ref{LemDerandom}}
\label{AppLemDerandom}

The main idea is to pass from probability statements to expectation
statements. We prove the forward direction, as the converse admits a
similar proof.

Consider the ``good'' event
\begin{align*}
\event = \eventwrap{W_{\nulls,(\lstar + 1)} \leq
  \transf\big(W_{\sigs,(\kstar - \lstar)}\big)},
\end{align*}
as well as the two ``bad'' events
\begin{align*}
\event_{0} & = \eventwrap{W_{\nulls,(\lstar + 1)} <
  \E\big[W_{\nulls,(\lstar + 1)}\big] - \ordD_{\nulls,\lstar +
    1}\big(\frac{\epsilon}{3}\big)}, \qquad \mbox{and} \\
\event_{1} & = \eventwrap{W_{\sigs,(\kstar - \lstar)} >
  \E\big[W_{\sigs,(\kstar - \lstar)}\big] + \ordD_{\sigs,\kstar -
    \lstar}\big(\frac{\epsilon}{3}\big)} .
\end{align*}
By the maximality of $\lstar$, we have $\P(\event) \geq \epsilon$,
while the definition of the concentration functions ensures
$\max\big\lbrace \P(\event_{0}),~\P(\event_{1})\big\rbrace \leq
\frac{\epsilon}{3}$.  Thus, if we define the event
\mbox{$\event_{\ast} = \event\setminus\big(\event_{0} \cup
  \event_{1}\big)$,} we are guaranteed that $\mprob[\event_\ast] \geq
\frac{\epsilon}{3}$.

Conditioned on $\event_{\ast}$, we have
\begin{align*}
\E\big[W_{\nulls,(\lstar + 1)}\big] - \ordD_{\nulls,\lstar +
  1}\big(\frac{\epsilon}{3}\big) \leq W_{\nulls,(\lstar + 1)} \leq
W_{\sigs,(\kstar - \lstar)} \leq \E\big[W_{\sigs,(\kstar -
    \lstar)}\big] + \ordD_{\sigs,\kstar -
  \lstar}\big(\frac{\epsilon}{3}\big).
\end{align*}
Comparing the left-hand side and right-hand side of this string of
inequalities yields the desired conclusion.


\subsection{Proof of Lemma~\ref{LemTransfer}}
\label{AppLemTransfer}

We first observe that the ``good'' event
\begin{align*}
\event & \defn \eventwrap{W_{\nulls,(\lstar)} > f(W_{\sigs,(\kstar -
    \lstar + 1)})}
\end{align*}
satisfies $\P(\event) \geq 1 - \frac{\epsilon}{3}$.  In order to
establish the claim, we need to show that the corresponding event for
the primed model, namely
\begin{align*}
\event' & = \eventwrap{W_{\nulls,(\lstar)}' > f(W_{\sigs,(\kstar -
    \lstar + 1)}) - \ordD_{0} - \ordD_{1}}
\end{align*}
satisfies $\P(\event') \geq 1 - \epsilon$. If this claim is proven,
the conclusion of the lemma will follow from the maximality of
${\lstar}'$ (cf. equation~\eqref{EqnFalseProxy}).

In order to establish the latter claim, we observe that $\event'
\supset \event \setminus \big(\event_{0} \cup \event_{1}\big)$, where
the bad events are defined by
\begin{align*}
\event_{0} & = \eventwrap{W_{\nulls,(\lstar)}' < W_{\nulls,(\lstar)} -
  \ordD_{0}}, \\ \event_{1} & = \eventwrap{W_{\sigs,(\kstar - \lstar +
    1)}' < W_{\sigs,(\kstar - \lstar + 1)} - \ordD_{1}} .
\end{align*}
Given these definitions, the inclusion is clear. Likewise, it is
immediate from the assumption of closeness that $\max\big\lbrace
\P(\event_{0}), \P(\event_{1})\big\rbrace \leq \frac{\epsilon}{3}$, so
that
\begin{align*}
\P\bigg(\event \setminus \big(\event_{0} \cup \event_{1}\big)\bigg) &
\geq 1 - \epsilon .
\end{align*}


\subsection{Proof of Lemma~\ref{LemDecouple}}
\label{AppLemDecouple}
    
The proof is similar to that of Lemmas~\ref{LemDerandom}
and~\ref{LemTransfer}. In this case, the ``good'' event is given by
\begin{align*}
\event'' = \eventwrap{W_{\nulls,({\lstar}')} >
  \transf(W_{\sigs,(\kstar - {\lstar}' + 1)}) - 4\Delta} .
\end{align*}
If we can show that $\P\big(\event''\big) \geq 1 - \epsilon$, the
conclusion of the lemma will follow from the maximality of ${\lstar}'$
and ${\lstar}''$ (see definition~\eqref{EqnFalseProxy}). Let
\begin{align*}
\event' = \eventwrap{W_{\nulls,({\lstar}')} > \transf(W_{\sigs,(\kstar - {\lstar}' + 1)}) }
\end{align*}
denote the corresponding event for the primed model $\model'$. Note
that by the definition of ${\lstar}'$, we have $\P\big(\event'\big)
\geq 1 - \frac{\epsilon}{3}$.

In order to control $\mprob[\event'']$, consider as usual the ``bad''
events
\begin{align*}
\event_{0} & = \eventwrap{\big|W_{\nulls,({\lstar}')} -
  W_{\sigs,({\lstar}')}'\big| \geq
  2\Delta_{\nulls,{\lstar}'}\big(\epsilon/6\big)}, \\ \event_{1} & =
\eventwrap{\big|W_{\sigs,(\kstar - {\lstar}' + 1)} - W_{\sigs,(\kstar
    - {\lstar}' + 1)}\big| \geq 2\Delta_{\sigs,\kstar - {\lstar}' +
    1}} .
\end{align*}
By two applications of Case I of the concentration
condition~\eqref{EqnConcentration}, we find that
\begin{align*}
\max\big\lbrace
\P\big(\event_{0}\big),~\P\big(\event_{1}\big)\big\rbrace \leq
\frac{\epsilon}{3}.
\end{align*}
Given the set inclusion $\event'' \supset
\event'\setminus\big(\event_{0} \cup \event_{1}\big)$, we conclude
that
\begin{align*}
\P\big(\event''\big) \geq 1 - \frac{\epsilon}{3} - \frac{2\epsilon}{3}
= 1 - \epsilon,
\end{align*}
as claimed.


\section{Details of simulations}
\label{SecSimulation}

In this section, we describe how to construct curves like the ones in
Figure~\ref{FigPredictedCurves}. For this, assume we have chosen a
fixed number $S$ of points at which to sample the curves. Given this
choice, we estimate both our theoretical lower bound (from
Theorem~\ref{ThmProxies}) on FNR and the attained FNR of BH at each of
the following FDR levels:
\begin{align*}
\Qs = \left \lbrace \frac{\epsilon}{B} \cdot b \; \mid \; b \in [0, B)
  \right \rbrace .
\end{align*}
We repeat this procedure for each model $\model$ in a set
$\allmodels$. As a result, we obtain predicted points on the
lower-left boundary of the feasible region, denoted as
\begin{align*}
\bigg((\fdr_{b},~\fnr_{\lo,b,\model})\bigg)_{b = 1, \dots, B, \model
  \in \allmodels}
\end{align*}
as well as points on the actual FDR-FNR tradeoff curve for BH, denoted
as
\begin{align*}
  \bigg((\fdr_{b},~ \fnr_{\BH,b,\model})\bigg)_{b = 1, \dots, B,
    \model \in \allmodels}.
\end{align*}
Based on these outputs, we estimate a value for a single
model-independent constant $c$ such that multiplying the $\fdr_{\lo}$
values by $c$ yields the best fit to the actual BH values, while still
preserving the lower-bounding property. The objective we use is a
simple least-squares objective given by
\begin{align*}
\mathcal{L}(c) = \sum_{\model \in \allmodels} \sum_{b = 1}^{B}
\bigg(\fdr_{\BH,b,\model} - \fdr_{\lo,b,\model} \cdot c\bigg)^{2} .
\end{align*}
Fortunately, the optimal estimate $\hat{c}$ subject to the
lower-bounding constraint can be found in closed form via
\begin{align*}
\hat{c}_{0} & = \frac{\sum_{b, \model}
  \fdr_{\BH,b,\model}}{\sum_{b,\model} \fdr_{\lo,b,\model}},
\\ \hat{c} & = \begin{cases} \hat{c}_{0} & ~\text{if}~\hat{c}_{0} \leq
  \min_{s,\model} \frac{\fdr_{\BH,b,\model}}{\fdr_{\lo,b,\model}},
  \\ \min_{s,\model} \frac{\fdr_{\BH,b,\model}}{\fdr_{\lo,b,\model}}
  &~\text{otherwise.}
\end{cases}
\end{align*}

For the experiment used to generate Figure~\ref{FigPredictedCurves},
we used $B = 25$, $\epsilon = 0.25$, and a set of sparse Gaussian
sequence models with locations shifts $\mu$ based on scalings $m =
n^{1 - \betapar}$ and $\mu = \sqrt{2r \log{n}}$. For all models, we
set $n = 10000$, while $\betapar$ varied within $\lbrace 0.5, 0.6, 0.7
\rbrace$ and $r$ varied within $\lbrace \betapar + 0.01,~\betapar +
0.05,~\betapar + 0.1\big\rbrace$ for each setting of $\betapar$.


\end{document}